\documentclass[10pt]{amsart}
\usepackage[T1]{fontenc}
\usepackage{lmodern}
\usepackage{geometry}
\usepackage[latin1] {inputenc}
\usepackage{amsmath}
\usepackage{amsfonts, amssymb, textcomp}
\usepackage[colorlinks=flase, linkcolor=red,urlcolor=green, citecolor=blue]{hyperref}
\usepackage{subeqnarray}
\usepackage[usenames]{color}
\usepackage{latexsym}
\usepackage{fancyhdr}
\usepackage{longtable}
\usepackage{amsmath, amssymb}
\usepackage{graphicx}

\numberwithin{equation}{section}

\theoremstyle{plain}
\newtheorem{theorem}{Theorem}[section]
\newtheorem{lemma}[theorem]{Lemma}
\newtheorem{proposition}[theorem]{Proposition}
\newtheorem{corollary}[theorem]{Corollary}
\theoremstyle{definition}
\newtheorem{definition}[theorem]{Definition}
\newtheorem{example}[theorem]{Example}
\newtheorem{remark}[theorem]{Remark}

\newcommand{\RR}{\mathbb{R}}
\newcommand{\CC}{\mathbb{C}}
\newcommand{\NN}{\mathbb{N}}
\newcommand{\R}{\mathbb{R}}
\newcommand{\C}{\mathbb{C}}
\newcommand{\N}{\mathbb{N}}

\let\on=\operatorname

\newenvironment{demo}[1]{\par\smallskip\noindent{\bf #1.}}{\par\smallskip}

\title[Sectorial extensions, via Laplace transforms, in ultraholomorphic classes]{Sectorial extensions, via Laplace transforms, in ultraholomorphic classes defined by weight functions}

\author[J.~Jim\'{e}nez-Garrido, J.~Sanz, and G.~Schindl]{Javier Jim\'{e}nez-Garrido, Javier Sanz and Gerhard Schindl}

\begin{document}
\begin{abstract}
We prove several extension theorems for Roumieu ultraholomorphic classes of functions in sectors of the Riemann surface of the logarithm which are defined by means of a weight function or weight matrix. Our main aim is to transfer the results of V. Thilliez
from the weight sequence case to these different, or more general, frameworks. The technique rests on the construction of suitable kernels for a truncated Laplace-like integral transform, which provides the solution without resorting to Whitney-type extension results for ultradifferentiable classes. As a byproduct, we obtain an extension in a mixed weight-sequence setting in which assumptions on the sequence are minimal.
\end{abstract}

\date{\today}

\maketitle

%
%
%

%
%

\section{Introduction}

The main aim of this paper is to prove the surjectivity of the Borel map (via the existence of right inverses for this map) in ultraholomorphic classes of functions in unbounded sectors defined by means of weight functions or weight matrices, so generalizing to this framework previous results available only in the ultradifferentiable setting. Let us start by reviewing such results and motivating our approach.

Ultradifferentiable classes of smooth functions in sets of $\mathbb{R}^n$, defined by suitably restricting the growth of their derivatives, have been extensively studied since the beginning of the 20th century. In particular, the study of the injectivity and surjectivity of, or the existence of right inverses for, the Borel map (respectively, the Whitney map), sending a function in this class to the family of its derivatives at a given point (resp., at every point in a given closed subset of $\mathbb{R}^n$), has attracted much attention. In case the restriction of growth is specified in terms of a sequence of positive real numbers, the corresponding classes are named after Denjoy and Carleman, who characterized the injectivity of the Borel map back in 1923.
The surjectivity of, and the existence of right inverses for, the Borel map was solved 1988 by H.-J. Petzsche~\cite{petzsche}, and the Whitney extension result was treated by J. Chaumat and A. M. Chollet~\cite{ChaumatChollet94}.
From the seminal work of R. W. Braun, R. Meise and B. A. Taylor~\cite{BraunMeiseTaylor90}, who modified the original approach of A. Beurling, it is also standard to consider classes in which the growth control is made by a weight function, whose properties allow one to conveniently apply Fourier analysis in this setting thanks to suitable Paley-Wiener-like results. The study of the surjectivity of the Borel and Whitney maps and their right inverses in this situation was done in the 1980's and 1990's by several authors, we mention J. Bonet, R. W. Braun, J. Bruna, M. Langenbruch, R. Meise and B. A. Taylor (see~\cite{MeiseTaylor89,BonetBraunMeiseTaylorWhitneyextension} and the references therein). A last step in this context has been recently taken by A. Rainer and the third author
\cite{dissertation,compositionpaper}, who considered classes defined by weight matrices, what strictly includes both the Denjoy-Carleman and the Braun-Meise-Taylor approaches, and also obtained results in the same line~\cite{dissertation,whitneyextensionweightmatrix}.

However, the study of similar problems for classes of holomorphic functions is much more recent, and it has been motivated by the increasing interest on asymptotic expansions, a theory put forward by H. Poincar\'e at the end of the 19th century. In order to give a full analytical meaning to the formal power series solutions of meromorphic linear systems of ordinary differential equations at an irregular singular point in the complex domain, in the 1980's J. P. Ramis, B. Malgrange, Y. Sibuya and W. Balser, just to name a few, refined this concept by considering Gevrey asymptotic expansions of order $\alpha>1$: Its existence  for a function $f$, holomorphic in a sector $S$ (with vertex at 0) of the Riemann surface of the logarithm, amounts to the estimations $|f^{(n)}(z)|\le CA_n(n!)^{\alpha}$ in proper subsectors of $S$, for suitable $C,A>0$. This fact makes evident the close link between ultradifferentiable classes and those similarly introduced for holomorphic functions defined in sectors, which are called ultraholomorphic classes. In asymptotic theory it is also important to decide about the injectivity or surjectivity of the Borel map, sending a function to the sequence of its derivatives at the vertex (defined by an obvious limiting process). While the injectivity for Gevrey classes was already studied by Watson and Nevanlinna in the 1920's, the surjectivity result, known as Borel-Ritt-Gevrey theorem, is due to B. Malgrange (see \cite{Ramis,RamisSeriesksommables}), and V. Thilliez \cite{ThilliezExtGevrey} obtained right inverses for the Borel map. Corresponding results for Gevrey functions in several variables were obtained by Y. Haraoka~\cite{Haraoka} and the second author~\cite{SanzGevreyseveralvar}. For general Denjoy-Carleman ultraholomorphic classes in unbounded sectors, in which the sequence $((n!)^{\alpha})_n$ is replaced by a general sequence $M=(M_n)_n$ subject to standard assumptions, the first results on the surjectivity of the Borel map and the existence of right inverses were obtained in 2000 by J. Schmets and M. Valdivia~\cite{Schmetsvaldivia00}, and these were improved in some respects by V. Thilliez~\cite{Thilliezdivision}. In this last paper, a growth index $\gamma(M)$ associated with the sequence $M$ plays a crucial role, limiting from above the opening of the sector for which extension operators are shown to
exist. Finally, in case the sequence $M$ admits a proximate order definitive results for injectivity and surjectivity were obtained by the second author in~\cite{Sanzflatultraholomorphic}, and a forthcoming paper~\cite{JimenezSanzSchindlInjSurj} will completely solve the injectivity problem for general logarithmically convex sequences, and it will provide significantly improved information for the surjectivity as long as strongly regular sequences are considered. However, no attempt has been made so far to study these problems for ultraholomorphic classes defined by weight functions or matrices, and our present paper is a first step in this direction.

The main ingredient for our construction of extension operators is the use of a truncated Laplace-like integral transform whose kernel is obtained from optimal flat functions, i.e., functions which are not only flat (in the sense that they have a null asymptotic expansion, and so an exponential decrease in terms of the sequence $M$) but admit also exponential estimates from below. This technique rests on the fundamental idea of B. Malgrange, and it has already been fruitful in an alternative proof by A. Lastra, S. Malek and the second author~\cite{LastraMalekSanzContinuousRightLaplace} of the extension results of V. Thilliez~\cite{Thilliezdivision}, and also in~\cite{Sanzflatultraholomorphic}. While the construction of sectorially (optimal) flat functions is an adaptation of the ideas by V. Thilliez, we will not use any Whitney-type extension result from the ultradifferentiable setting: A suitable integral kernel is defined from the flat functions available, and its moments are proved to be estimated from above and below by sequences belonging to the weight matrix defining the ultraholomorphic class (see Proposition~\ref{propmequivm}). The opening of the sectors for which the construction is possible is again controlled by a new growth index $\gamma(\omega)$, associated in this case with the defining weight function $\omega$, and which allows one to turn qualitative properties of $\omega$ into quantitative ones (see in this respect the Lemmas~\ref{lemmagammaomegapositiveimpliesomega1} and~\ref{strongweightspace}). For a detailed information about this and other indices for $\omega$, and their relation to the indices $\gamma(M)$ of Thilliez or $\omega(M)$ (introduced in~\cite{Sanzflatultraholomorphic}), we refer to a paper in preparation~\cite{JimenezSanzSchindlIndices}.

The main result, Theorem~\ref{comparisonmaintheorem}, states the surjectivity of the Borel map in ultraholomorphic classes defined by a weight function $\omega$, satisfying some standard assumptions and such that $\gamma(\omega)>1$, as long as the opening of the sector is smaller than $\pi(\gamma(\omega)-1)$.
As an intermediate step, Theorem~\ref{theoExtensionOperatorsMatrix} proves the surjectivity of the Borel map in ultraholomorphic classes associated with a weight matrix which is, in turn, obtained from a suitable weight function $\tau$ with $\gamma(\tau)>0$, and in sectors of opening smaller than $\pi\gamma(\tau)$. The close link between these two statements is given by the upper and lower Legendre transforms, classical tools in the study of the Laplace and Fourier transforms and which (loosely speaking) allow one to move from $\omega$ to $\tau$ and viceversa.
Observe that, as a byproduct of the existence of optimal flat functions, we deduce that the Borel map is not injective for these classes in such narrow sectors.

It is worthy to mention that a different proof of Theorem 7.7 has been obtained by the authors, see \cite{sectorialextension}, by closely following the arguments in the proof of V. Thilliez for Carleman ultraholomorphic classes defined by strongly regular weight sequences: Whitney extension results from the ultradifferentiable setting are needed, and moreover a careful adaptation of a ramification process, also applied by Thilliez, has to be carried out for classes associated with weight matrices.
A distinctive feature of this new procedure is the possibility of adapting it in order to treat also the Beurling case, which has not been considered in the present paper.

A last paragraph in the paper is devoted to the implications of our main result when Denjoy-Carleman ultraholomorphic classes are considered. If the weight sequence $M$ is strongly regular we recover the result of Thilliez, but if we drop the moderate growth condition we are able to prove an extension result in a mixed setting, meaning that the weight sequence defining the class of sequences we depart from has to be changed into a precise, larger (nonequivalent one) weight sequence defining the ultraholomorphic class where the interpolating function dwells. It is worthy to emphasize that,
as illustrated by Example~\ref{examSequenceBadGrowth}, there do exist sequences which do not satisfy any of the standard growth properties assumed in previous extension results in the literature, but to which our results can be applied, see Subsection~\ref{subsectMixedSetting}.

The paper is organized as follows. Section~\ref{sectBasicDefs} contains all the preliminary, mostly well-known, information concerning weight sequences, weight functions and weight matrices, and it introduces the ultraholomorphic classes we will consider, among which those associated with weight functions or matrices are new in the literature. It ends with Lemma~\ref{lemmaexponentialdecrease.hm}, which will be important for rephrasing flatness in our ultraholomorphic classes by means of some standard auxiliary functions. In Section~\ref{sectLegendreConjWeightMatrtoWeightFunc} we recall some basic facts about Legendre (also called Young) conjugates
and obtain the equivalence of different weight functions constructed via these conjugates.
The information about Thilliez's growth index for a weight sequence, and about a new growth index for weight functions, is described in Section~\ref{growthindices}, and the next section is devoted to
the equality of different ultraholomorphic classes  associated with weight functions or matrices.
After a characterization of flat functions (Lemma~\ref{lemmaflatnessultraholclasses}), the construction of optimal flat functions is the aim of Section~\ref{sectSectoriallyFlatFunct}. Finally, Section~\ref{sectRightInver1var} is devoted to
the main results, Theorems~\ref{comparisonmaintheorem} and \ref{theoExtensionOperatorsMatrix}, a rephrasing of the latter
(Corollary~\ref{coroExtensionOperatorsWeights})
in terms of classes defined by weight functions, and a closing subsection about a mixed setting extension procedure for classes defined by weight sequences.

\section{Basic definitions}\label{sectBasicDefs}
This section is devoted to fixing some notations, introducing the main properties of weight sequences, functions or matrices which we will deal with, and defining the ultraholomorphic classes of Roumieu type under consideration.

We denote by 
$\mathcal{H}$ the class of holomorphic functions. We will write $\NN_{>0}=\{1,2,\dots\}$ and $\NN=\NN_{>0}\cup\{0\}$, moreover we put $\RR_{>0}:=\{x\in\RR: x>0\}$, i.e. the set of all positive real numbers.


\subsection{Weight sequences}\label{subsectionWeightSequences}



A sequence $M=(M_k)_k\in\RR_{>0}^{\NN}$ is called a {\itshape weight sequence}. We define also $m=(m_k)_k$ by
$$m_k:=\frac{M_k}{k!},\quad k\in\NN,
$$
and $\mu=(\mu_k)_k$ by
$$\mu_0:=1;\quad\mu_k:=\frac{M_k}{M_{k-1}},\ \ k\in\NN_{>0}.
$$
%

$M$ is called {\itshape normalized} if $1=M_0\le M_1$ (this condition may always be assumed without loss of generality).

We list now some interesting and standard properties for weight sequences:

$(1)$ $M$ is {\itshape log-convex}, if
$$\hypertarget{lc}{(\text{lc})}:\Leftrightarrow\;\forall\;j\in\NN_{>0}:\;M_j^2\le M_{j-1} M_{j+1}$$
and {\itshape strongly log-convex}, if
$$\hypertarget{slc}{(\text{slc})}:\Leftrightarrow\;\forall\;j\in\NN_{>0}:\;m_j^2\le m_{j-1} m_{j+1}.$$
We recall that for every weight sequence $M=(M_k)_k\in\RR_{>0}^{\NN}$ one has
$$
\liminf_{k\to\infty}\mu_k\le\liminf_{k\to\infty}(M_k)^{1/k}\le
\limsup_{k\to\infty}(M_k)^{1/k}\le\limsup_{k\to\infty}\mu_k.$$
If $M$ is log-convex and normalized, then $M$, $((M_k)^{1/k})_{k\in\NN}$ and $(\mu_k)_{k\in\NN}$ are nondecreasing, and so
$\lim_{k\rightarrow\infty}(M_k)^{1/k}=+\infty$ if, and only if, $\lim_{k\rightarrow\infty}\mu_k=+\infty$.
Moreover, $M_jM_k\le M_{j+k}$ holds for all $j,k\in\NN$, e.g. see \cite[Remark 2.0.3, Lemmata 2.0.4, 2.0.6]{diploma}.

$(2)$ $M$ has {\itshape moderate growth} if
$$\hypertarget{mg}{(\text{mg})}:\Leftrightarrow\exists\;C\ge 1\;\forall\;j,k\in\NN:\;M_{j+k}\le C^{j+k}M_j M_k.$$
Note that, by elementary estimates, $M$ has \hyperlink{mg}{(\text{mg})} if, and only if, $m$ has \hyperlink{mg}{(\text{mg})}.



$(3)$ $M$ has
$(\gamma_1)$ if
$$\hypertarget{gamma1}{(\gamma_1)}:\Leftrightarrow\sup_{p\in\NN_{>0}}\frac{\mu_p}{p}\sum_{k\ge p}\frac{1}{\mu_k}<+\infty.$$
In the literature \hyperlink{gamma1}{$(\gamma_1)$} is also called ``strong nonquasianalyticity condition''. 

Due to technical reasons it is often convenient to assume several properties for $M$ at the same time and hence we define the class

\centerline{$M\in\hypertarget{SRset}{\mathcal{SR}}$, if $M$ is normalized and has \hyperlink{slc}{$(\text{slc})$}, \hyperlink{mg}{$(\text{mg})$} and \hyperlink{gamma1}{$(\gamma_1)$}.}

Using this notation we see that $M\in\hyperlink{SRset}{\mathcal{SR}}$ if and only if $m$ is a {\itshape strongly regular sequence} in the sense of \cite[1.1]{Thilliezdivision} (and this terminology has also been used by several authors so far, e.g. see \cite{Sanzflatultraholomorphic}, \cite{Sanzsummability}). At this point we want to make the reader aware that in \cite{Thilliezdivision} a slightly different notation and terminology is used, due to the fact that the main role in the statements there is assigned to the sequence which here is denoted by $m$, and not to the sequence denoted here by $M$.




%
%

We write $M\le N$, and we say $M$ is a minorant of $N$, and $N$ a majorant of $M$, if and only if $M_p\le N_p$ holds for all $p\in\NN$, and define
$$M\hypertarget{precsim}{\precsim}N:\Leftrightarrow\;\exists \;C\ge 1\;\forall\;p\in\NN:\; M_p\le C^pN_p\Longleftrightarrow\sup_{p\in\NN_{>0}} \left(\frac{M_p}{N_p}\right)^{1/p}<+\infty.$$
$M$ and $N$ are called {\itshape equivalent} if
$$M\hypertarget{approx}{\approx}N:\Leftrightarrow\; M\hyperlink{precsim}{\precsim}N\;\text{and}\;
N\hyperlink{precsim}{\precsim}M.$$
Moreover, if we write $\nu_0:=1$, $\nu_p:=N_p/N_{p-1}$, $p\in\NN_{>0}$, we introduce the stronger relation
$$M\hypertarget{mupreceq}{\preceq}N:\Leftrightarrow\;\exists\;C\ge 1\;\forall\;p\in\NN:\;\mu_p\le C\nu_p\Longleftrightarrow\sup_{p\in\NN}\frac{\mu_p}{\nu_p}<+\infty$$
and call them {\itshape strongly equivalent} if
$$M\hypertarget{simeq}{\simeq}N:\Leftrightarrow\;
M\hyperlink{mupreceq}{\preceq}N\;\text{and}\;
N\hyperlink{mupreceq}{\preceq}M.$$
If we write $n=(n_k)_k$ for
$n_k:=\frac{N_k}{k!}$, $k\in\NN$, then it is clear that
$M\hyperlink{precsim}{\precsim}N$ if, and only if, $m\hyperlink{precsim}{\precsim}n$, and that $M\hyperlink{mupreceq}{\preceq}N$ if, and only if, $m\hyperlink{mupreceq}{\preceq}n$.

Define the set
$$\hypertarget{LCset}{\mathcal{LC}}:=\{M\in\RR_{>0}^{\NN}:\;M\;\text{normalized, log-convex},\;\lim_{k\rightarrow\infty}(M_k)^{1/k}=+\infty\}.$$
We warn the reader that in previous works by the authors~\cite{JimenezGarridoSanz,logconvexnonproximate} the condition $\lim_{k\rightarrow\infty}(M_k)^{1/k}=+\infty$ was equivalently expressed (see above) as $\lim_{k\rightarrow\infty}\mu_k=+\infty$.

The {\itshape Gevrey sequence} of order $s\ge 1$ will be denoted by $G^s:=(p!^s)_p$, for $s>1$ it satisfies all properties listed above.

\subsection{Weight functions $\omega$ in the sense of Braun-Meise-Taylor}\label{weightfunctionclasses}
A function $\omega:[0,\infty)\rightarrow[0,\infty)$ 
is called a \textit{weight function} if it is continuous, nondecreasing, $\omega(0)=0$ and $\lim_{x\rightarrow\infty}\omega(x)=+\infty$.

In case we also have $\omega(x)=0$ for all $x\in[0,1]$, we say $\omega$ is a \textit{normalized weight}.


Moreover we consider the following conditions:
\begin{itemize}
\item[\hypertarget{om1}{$(\omega_1)}$] $\omega(2t)=O(\omega(t))$ as $t\rightarrow+\infty$.


\item[\hypertarget{om3}{$(\omega_3)$}] $\log(t)=o(\omega(t))$ as $t\rightarrow+\infty$ ($\Leftrightarrow\lim_{t\rightarrow+\infty}\frac{t}{\varphi_{\omega}(t)}=0$).

\item[\hypertarget{om4}{$(\omega_4)$}] The function $\varphi_{\omega}:\RR\to\RR$, given by $\varphi_{\omega}(t)=\omega(e^t)$, is a convex function on $\RR$.

\item[\hypertarget{om5}{$(\omega_5)$}] $\omega(t)=o(t)$ as $t\rightarrow+\infty$.

\item[\hypertarget{om6}{$(\omega_6)$}] $\exists\;H\ge 1\;\forall\;t\ge 0:\;2\omega(t)\le\omega(H t)+H$.



\item[\hypertarget{omsnq}{$(\omega_{\text{snq}})$}] $\exists\;C>0:\;\forall\;y>0: \int_1^{\infty}\frac{\omega(y t)}{t^2}dt\le C\omega(y)+C$.

\end{itemize}

We mention that this list of properties is extracted from a larger one used already in~\cite{dissertation}, what explains the lack of $(\omega_2)$, irrelevant in this paper.

An interesting example is the weight function $\sigma_s(t):=\max\{0,\log(t)^s\}$, $s>1$, which satisfies all listed properties except \hyperlink{om6}{$(\omega_6)$}. It is well-known that the weight $t\mapsto t^{1/s}$ yields the Gevrey class $G^s$ of index $s>1$, it satisfies all listed properties (except normalization).

%



For a normalized weight $\omega$ satisfying \hyperlink{om3}{$(\omega_3)$}
we define the {\itshape Legendre-Fenchel-Young-conjugate}
\begin{equation}\label{legendreconjugate}
\varphi^{*}_{\omega}(x):=\sup\{x y-\varphi_{\omega}(y): y\in\R\}=\sup\{x y-\varphi_{\omega}(y): y\ge 0\},\;\;\;x\ge 0,
\end{equation}
with the following properties, e.g. see \cite[Remark 1.3, Lemma 1.5]{BraunMeiseTaylor90}: It is nonnegative, convex and nondecreasing, $\varphi^{*}_{\omega}(0)=0$, the map $x\mapsto\frac{\varphi^{*}_{\omega}(x)}{x}$ is nondecreasing in $[0,+\infty)$ and $\lim_{x\rightarrow\infty}\frac{\varphi^{*}_{\omega}(x)}{x}= \infty$.

Moreover, $\omega$ has also \hyperlink{om4}{$(\omega_4)$} if, and only if, $\varphi^{**}_{\omega}=\varphi_{\omega}$,  and then the map $x\mapsto\frac{\varphi_{\omega}(x)}{x}$ is also nondecreasing in $[0,+\infty)$.\par

\begin{remark}\label{remarkConditionomega4}
It is interesting to note, as it was done in~\cite[p.\ 15]{dissertation}, that condition \hyperlink{om4}{$(\omega_4)$}, appearing in~\cite{BraunMeiseTaylor90}, was necessary in order to show that certain classes of compactly supported functions defined by decay properties of their Fourier transform in terms of a weight function $\omega$ could be alternatively represented as those consisting of functions whose derivatives' growth may be controlled by the
Legendre-Fenchel-Young-conjugate of $\omega$. Since we will work in a different framework, we will only assume this condition whenever the equality $\varphi^{**}_{\omega}=\varphi_{\omega}$ is needed in our arguments.
\end{remark}

Given a weight function $\omega$ and $s>0$, we define a new weight function $\omega^s$ by
\begin{equation}\label{defiRamifiedWeightFunction}
\omega^s(t):=\omega(t^s),\quad t\ge 0.
 \end{equation}
Also, given a function $\tau:(0,\infty)\to\mathbb{R}$ we sometimes consider the function $\tau^{\iota}$ given by $\tau^{\iota}(t):=\tau(1/t)$, $t>0$.

If $\omega$ satisfies any of the properties \hyperlink{om1}{$(\omega_1)$}, \hyperlink{om3}{$(\omega_3)$}, \hyperlink{om4}{$(\omega_4)$} or \hyperlink{om6}{$(\omega_6)$}, then the same holds for $\omega^s$, but \hyperlink{om5}{$(\omega_5)$} or \hyperlink{omsnq}{$(\omega_{\text{snq}})$} might not be preserved. Indeed, this last fact motivates the introduction of the index $\gamma(\omega)$ in this paper, see Subsection~\ref{growthindexgamma}.

Let $\sigma,\tau$ be weight functions, 
we write
$$\sigma\hypertarget{ompreceq}{\preceq}\tau:\Leftrightarrow\tau(t)=O(\sigma(t))\;\text{as}\;t\rightarrow+\infty$$
and call them \textit{equivalent} if
$$\sigma\hypertarget{sim}{\sim}\tau:\Leftrightarrow\sigma\hyperlink{ompreceq}{\preceq}\tau\;\text{and}\;\tau\hyperlink{ompreceq}{\preceq}\sigma.$$

%

We recall \cite[Proposition 1.3]{MeiseTaylor88}, where \hyperlink{omsnq}{$(\omega_{\text{snq}})$} was characterized, and \cite[Corollary 1.4]{MeiseTaylor88}:

\begin{proposition}\label{Prop13MT88}
Let $\omega:[0,+\infty)\longrightarrow[0,+\infty)$ be a weight function.
The following are equivalent::
\begin{itemize}
\item[$(i)$] $\lim_{\varepsilon\rightarrow 0}\limsup_{t\rightarrow+\infty}\frac{\varepsilon\omega(t)}{\omega(\varepsilon t)}=0$,
\item[$(ii)$] $\exists\;K>1$ such that $\limsup_{t\rightarrow+\infty}\frac{\omega(Kt)}{\omega(t)}<K$,
\item[$(iii)$] $\omega$ satisfies \hyperlink{omsnq}{$(\omega_{\on{snq}})$},
\item[$(iv)$] There exists a nondecreasing concave function $\kappa: [0,+\infty)\longrightarrow[0,+\infty)$ such that $\omega\hyperlink{sim}{\sim}\kappa$ and $\kappa$ satisfies \hyperlink{omsnq}{$(\omega_{\on{snq}})$}. More precisely $\kappa=\kappa_{\omega}$ with
    \begin{equation*}
\kappa_{\omega}(t):=\int_1^{\infty}\frac{\omega(tu)}{u^2}du=t\int_t^{\infty}\frac{\omega(u)}{u^2}du,\;\;\;\forall\;t>0\hspace{20pt}\kappa_{\omega}(0)=0.
\end{equation*}
\end{itemize}
Consequently, $\omega$ has also \hyperlink{om1}{$(\omega_1)$} and \hyperlink{om5}{$(\omega_5)$}. If $\omega$ satisfies one of the equivalent conditions above, then there exists some $0<\alpha<1$ such that $\omega(t)=O(t^{\alpha})$ as $t\rightarrow\infty$.
\end{proposition}

It is well-known that each of the properties \hyperlink{om1}{$(\omega_1)$}, \hyperlink{om3}{$(\omega_3)$} or \hyperlink{om4}{$(\omega_4)$}
can be transferred from $\omega$ to $\kappa_{\omega}$, see e.g. \cite[Remark 3.2]{BonetMeiseTaylorSurjectivity}.

Note that concavity of a weight function $\omega$ implies subadditivity
(i.e. $\omega(s+t)\le \omega(s)+\omega(t)$ for every $s,t\ge 0$; the proof needs the fact that $\omega(0)=0$), and this in turn yields \hyperlink{om1}{$(\omega_1)$}.



\subsection{
Weight matrices}\label{classesweightmatrices}
For the following definitions and conditions see also \cite[Section 4]{compositionpaper}.

Let $\mathcal{I}=\RR_{>0}$ denote the index set, a {\itshape weight matrix} $\mathcal{M}$ associated with $\mathcal{I}$ is a (one parameter) family of weight sequences $\mathcal{M}:=\{M^x\in\RR_{>0}^{\NN}: x\in\mathcal{I}\}$, such that
$$\hypertarget{Marb}{(\mathcal{M})}:\Leftrightarrow\;\forall\;x\in\mathcal{I}:\;M^x\;\text{is normalized, nondecreasing},\;M^{x}\le M^{y}\;\text{for}\;x\le y.$$
We call a weight matrix $\mathcal{M}$ {\itshape standard log-convex,} if
$$\hypertarget{Msc}{(\mathcal{M}_{\on{sc}})}:\Leftrightarrow(\mathcal{M})\; \text{and}\;\forall\;x\in\mathcal{I}:\;M^x\in\hyperlink{LCset}{\mathcal{LC}}.$$
Moreover, we put $m^x_p:=\frac{M^x_p}{p!}$ for $p\in\NN$, and $\mu^x_p:=\frac{M^x_p}{M^x_{p-1}}$ for $p\in\NN_{>0}$, $\mu^x_0:=1$.

A matrix is called {\itshape constant} if $\mathcal{M}=\{M\}$ or more generally if $M^x\hyperlink{approx}{\approx}M^y$ for all $x,y\in\mathcal{I}$.

We are going to consider the following properties for $\mathcal{M}$:\par\vskip.3cm

\hypertarget{R-mg}{$(\mathcal{M}_{\{\text{mg}\}})$} \hskip1cm $\forall\;x\in\mathcal{I}\;\exists\;C>0\;\exists\;y\in\mathcal{I}\;\forall\;j,k\in\NN: M^x_{j+k}\le C^{j+k} M^y_j M^y_k$.\par\vskip.3cm

\hypertarget{R-L}{$(\mathcal{M}_{\{\text{L}\}})$} \hskip1cm $\forall\;C>0\;\forall\;x\in\mathcal{I}\;\exists\;D>0\;\exists\;y\in\mathcal{I}\;\forall\;k\in\NN: C^k M^x_k\le D M^y_k$.\vskip.3cm

Let $\mathcal{M}=\{M^x: x\in\mathcal{I}\}$ and $\mathcal{N}=\{N^x: x\in\mathcal{J}\}$ be \hyperlink{Marb}{$(\mathcal{M})$}, define
$$\mathcal{M}\hypertarget{Mroumprecsim}{\{\precsim\}}\mathcal{N}\;:
\Leftrightarrow\forall\;x\in\mathcal{I}\;\exists\;y\in\mathcal{J}:\;
M^x\hyperlink{precsim}{\precsim}N^y,$$
and \textit{equivalence} of matrices,
$$\mathcal{M}\hypertarget{Mroumapprox}{\{\approx\}}\mathcal{N}\;:
\Leftrightarrow\mathcal{M}\hyperlink{Mroumprecsim}{\{\precsim\}}\mathcal{N}\;\text{and}\;\mathcal{N}\hyperlink{Mroumprecsim}{\{\precsim\}}\mathcal{M}$$


\subsection{Weight matrices obtained from weight functions}\label{weightmatrixfromfunction}
We summarize some facts which are shown in \cite[Section 5]{compositionpaper} and will be needed.
\begin{itemize}
\item[$(i)$] A central idea was that with each normalized weight function $\omega$ that has 
 $\hyperlink{om3}{(\omega_3)}$ 
    we can associate a \hyperlink{Msc}{$(\mathcal{M}_{\text{sc}})$} weight matrix $\Omega:=\{W^l=(W^l_j)_{j\in\NN}: l>0\}$ by

$$W^l_j:=\exp\left(\frac{1}{l}\varphi^{*}_{\omega}(lj)\right),
$$
    which moreover satisfies \hyperlink{R-mg}{$(\mathcal{M}_{\{\text{mg}\}})$}, more precisely
    \begin{equation}\label{newmoderategrowth}
    \forall\;l>0\;\forall\;j,k\in\NN:\;\;\;W^l_{j+k}\le W^{2l}_jW^{2l}_k.
    \end{equation}
    In general it is not clear that $W^x$ is strongly log-convex, i.e. $w^x$ is log-convex, too.

\item[$(ii)$] If $\omega$ has moreover $\hyperlink{om1}{(\omega_1)}$, then $\Omega$
    satisfies also
\hyperlink{R-L}{$(\mathcal{M}_{\{\text{L}\}})$}, more precisely
     \begin{equation}\label{newexpabsorb}
     \forall\;h\ge 1\;\exists\;A\ge 1\;\forall\;l>0\;\exists\;D\ge 1\;\forall\;j\in\NN:\;\;\;h^jW^l_j\le D W^{Al}_j.
     \end{equation}
     In fact we can take $A=(L(L+1))^a$, where $L\ge 1$ is the constant arising in \hyperlink{om1}{$(\omega_1)$}, i.e. $\omega(2t)\le L(\omega(t)+1)$, and $a\in\NN_{>0}$ is chosen minimal to have $\exp(a)\ge h$ (see the proof of \cite[Proposition 3.3.1, p.18]{dissertation} and of \cite[Lemma 5.9 $(5.10)$]{compositionpaper}).

\item[$(iii)$] Equivalent weight functions $\omega$ yield equivalent weight matrices with respect to  \hyperlink{Mroumapprox}{$\{\approx\}$}. Note that \hyperlink{R-mg}{$(\mathcal{M}_{\{\text{mg}\}})$} is stable with respect to \hyperlink{Mroumapprox}{$\{\approx\}$}, whereas \hyperlink{R-L}{$(\mathcal{M}_{\{\text{L}\}})$} not.

\item[$(iv)$] \hyperlink{om5}{$(\omega_5)$} implies $\lim_{p\rightarrow\infty}(w^l_p)^{1/p}=+\infty$ for all $l>0$. Conversely, $\lim_{p\rightarrow\infty}(w^l_p)^{1/p}=+\infty$ for all $l>0$ together with \hyperlink{om4}{$(\omega_4)$} imply \hyperlink{om5}{$(\omega_5)$}.

\end{itemize}


\begin{remark}\label{importantremark}
As can be seen in~\cite[Lemma\ 5.1.3]{dissertation}, if $\omega$ is a normalized weight function satisfying \hyperlink{om3}{$(\omega_3)$} and \hyperlink{om4}{$(\omega_4)$}, then $\omega$ satisfies
 \hyperlink{om6}{$(\omega_6)$} if, and only if, some/each $W^l$ satisfies \hyperlink{mg}{$(\on{mg})$}, and this amounts to the fact that $W^l\hyperlink{approx}{\approx}W^s$ for each $l,s>0$. Consequently, \hyperlink{om6}{$(\omega_6)$} is characterizing the situation when $\Omega$ is constant, i.e. all the weight sequences it consists of are equivalent to each other.
\end{remark}


\subsection{Classes of ultraholomorphic functions of Roumieu type}\label{ultraholomorphicroumieu}
For the following definitions, notation and more details we refer to \cite[Section 2]{Sanzflatultraholomorphic}. Let $\mathcal{R}$ be the Riemann surface of the logarithm. We wish to work in general unbounded sectors in $\mathcal{R}$ with vertex at 0, but all our results will be unchanged under rotation, so we will only consider sectors bisected by direction 0: For $\gamma>0$ we set $$S_{\gamma}:=\{z\in\mathcal{R}: |\arg(z)|<\frac{\gamma\pi}{2}\},
$$
i.e. the unbounded sector of opening $\gamma\pi$, bisected by direction $0$.

%

Let $M$ be a weight sequence, $S\subseteq\mathcal{R}$ an (unbounded) sector and $h>0$. We define
$$\mathcal{A}_{M,h}(S):=\{f\in\mathcal{H}(S): \|f\|_{M,h}:=\sup_{z\in S, p\in\NN}\frac{|f^{(p)}(z)|}{h^p M_p}<+\infty\}.$$
$(\mathcal{A}_{M,h}(S),\|\cdot\|_{M,h})$ is a Banach space and we put
\begin{equation*}
\mathcal{A}_{\{M\}}(S):=\bigcup_{h>0}\mathcal{A}_{M,h}(S).
\end{equation*}
$\mathcal{A}_{\{M\}}(S)$ is called the Denjoy-Carleman ultraholomorphic class (of Roumieu type) associated with $M$ in the sector $S$ (it is a $(LB)$ space). Analogously, we introduce the space of complex sequences
$$\Lambda_{M,h}:=\{a=(a_p)_p\in\CC^{\NN}: |a|_{M,h}:=\sup_{p\in\NN^n}\frac{|a_p|}{h^{p} M_{p}}<+\infty\}$$
and put $\Lambda_{\{M\}}:=\bigcup_{h>0}\Lambda_{M,h}$. The (asymptotic) {\itshape Borel map} $\mathcal{B}$ is given by
\begin{equation*}
\mathcal{B}:\mathcal{A}_{\{M\}}(S)\longrightarrow\Lambda_{\{M\}},\hspace{15pt}f\mapsto(f^{(p)}(0))_{p\in\NN},
\end{equation*}
where $f^{(p)}(0):=\lim_{z\in S, z\rightarrow 0}f^{(p)}(z)$.



Similarly as for the ultradifferentiable case, we now define ultraholomorphic classes associated with a normalized weight function $\omega$ satisfying $\hyperlink{om3}{(\omega_3)}$. Given an unbounded sector $S$, and for every $l>0$, we first define
$$\mathcal{A}_{\omega,l}(S):=\{f\in\mathcal{H}(S): \|f\|_{\omega,l}:=\sup_{z\in S, p\in\NN}\frac{|f^{(p)}(z)|}{\exp(\frac{1}{l}\varphi^{*}_{\omega}(lp))}<+\infty\}.$$
$(\mathcal{A}_{\omega,l}(S),\|\cdot\|_{\omega,l})$ is a Banach space and we put
\begin{equation*}
\mathcal{A}_{\{\omega\}}(S):=\bigcup_{l>0}\mathcal{A}_{\omega,l}(S).
\end{equation*}
$\mathcal{A}_{\{\omega\}}(S)$ is called the Denjoy-Carleman ultraholomorphic class (of Roumieu type) associated with $\omega$ in the sector $S$ (it is a $(LB)$ space). Correspondingly, we introduce the space of complex sequences
$$\Lambda_{\omega,l}:=\{a=(a_p)_p\in\CC^{\NN}: |a|_{\omega,l}:=\sup_{p\in\NN}\frac{|a_p|}{\exp(\frac{1}{l}\varphi^{*}_{\omega}(lp))}<+\infty\}$$
and put $\Lambda_{\{\omega\}}:=\bigcup_{l>0}\Lambda_{\omega,l}$. So in this case we get the Borel map $\mathcal{B}:\mathcal{A}_{\{\omega\}}(S)\longrightarrow\Lambda_{\{\omega\}}$.


Finally, we recall that ultradifferentiable function classes $\mathcal{E}_{\{\mathcal{M}\}}$, of {\itshape Roumieu type} and defined by a weight matrix $\mathcal{M}$, were introduced in \cite{dissertation}, see also \cite[4.2]{compositionpaper}.
Similarly, given a weight matrix $\mathcal{M}=\{M^x\in\RR_{>0}^{\NN}: x\in\RR_{>0}\}$ and a sector $S$ we may define ultraholomorphic classes  $\mathcal{A}_{\{\mathcal{M}\}}(S)$  of {\itshape Roumieu type} as
\begin{equation*}
\mathcal{A}_{\{\mathcal{M}\}}(S):=\bigcup_{x\in\RR_{>0}}\mathcal{A}_{\{M^x\}}(S),
\end{equation*}
and accordingly, $\Lambda_{\{\mathcal{M}\}}:=\bigcup_{x\in\RR_{>0}}\Lambda_{\{M^x\}}$.

It is straightforward to check that equivalent weight sequences, functions or matrices yield the same corresponding ultraholomorphic classes.

As said before in Subsection \ref{weightmatrixfromfunction}, if $\omega$ is a normalized weight function with $\hyperlink{om1}{(\omega_1)}$ and  $\hyperlink{om3}{(\omega_3)}$, the \hyperlink{Msc}{$(\mathcal{M}_{\text{sc}})$} weight matrix $\Omega:=\{W^l=(W^l_j)_{j\in\NN}: l>0\}$ given by $W^l_j:=\exp\left(\frac{1}{l}\varphi^{*}_{\omega}(lj)\right)$ satisfies \hyperlink{R-mg}{$(\mathcal{M}_{\{\text{mg}\}})$} (see~\eqref{newmoderategrowth}) and \hyperlink{R-L}{$(\mathcal{M}_{\{\text{L}\}})$} (see~\eqref{newexpabsorb}), and moreover \begin{equation}\label{equaEqualitySpacesWeightFunctionMatrix}
\mathcal{A}_{\{\omega\}}(S)=\mathcal{A}_{\{\Omega\}}(S)
\end{equation}
holds as locally convex vector spaces (this equality is an easy consequence of the way the seminorms are defined in these spaces and of property \hyperlink{R-L}{$(\mathcal{M}_{\{\text{L}\}})$}). As one also has $\Lambda_{\{\omega\}}=\Lambda_{\{\Omega\}}$, the Borel map $\mathcal{B}$ makes sense in these last classes,
$\mathcal{B}:\mathcal{A}_{\{\Omega\}}(S)\longrightarrow \Lambda_{\{\Omega\}}$.

In any of the considered ultraholomorphic classes, an element $f$ is said to be \textit{flat} if $f^{(p)}(0)=0$ for every $p\in\NN$, that is, $\mathcal{B}(f)$ is the null sequence.

With all this information we may claim: In general it seems reasonable to transfer well-known results and proofs from the weight sequence to the weight function setting, by using its associated weight matrix, and even to the setting in which an abstract weight matrix is given from the start.
Unfortunately in many cases it is often impossible to replace $M$ by $W^x$ in the proofs directly since, due to technical reasons, undesirable resp. too strong conditions have been imposed on the weight sequences.

\subsection{Functions $\omega_M$ and $h_M$}\label{assofunction}

Many properties of weight sequences may be easily studied by means of the auxiliary functions which we introduce now.

Let $M\in\RR_{>0}^{\NN}$ ($M_0=1$), then the {\itshape associated function} $\omega_M: [0,\infty)\rightarrow\RR\cup\{+\infty\}$ is defined by
\begin{equation*}
\omega_M(t):=\sup_{p\in\NN}\log\left(\frac{t^p}{M_p}\right)\;\;\;\text{for}\;t>0,\hspace{30pt}\omega_M(0):=0.
\end{equation*}
For an abstract introduction of the associated function we refer to \cite[Chapitre I]{mandelbrojtbook}, see also \cite[Definition 3.1]{Komatsu73}. If $\liminf_{p\rightarrow\infty}(M_p)^{1/p}>0$, then $\omega_M(t)=0$ for sufficiently small $t$, since $\log\left(\frac{t^p}{M_p}\right)<0\Leftrightarrow t<(M_p)^{1/p}$ holds for all $p\in\NN_{>0}$.

%

A basic assumption is $\lim_{p\rightarrow\infty}(M_p)^{1/p}=+\infty$, which implies that $\omega_M(t)<+\infty$ for any $t>0$, and so
$\omega_M$ is a weight function. If moreover $M$ is normalized, then $\omega_M$ also is.

Let $M\in\RR_{>0}^{\NN}$ with $\lim_{p\rightarrow\infty}(M_p)^{1/p}=+\infty$, then the sequence defined by
\begin{equation}\label{equaLCminorant}
M_p^{\on{lc}}:=\sup_{t>0}\frac{t^p}{\exp(\omega_M(t))},
\end{equation}
is the \textit{log-convex regularization}, or \textit{log-convex minorant}, of the sequence $M$ (i.e., it is the greatest among the log-convex minorants of the sequence $M$), see for example~\cite[(3.2)]{Komatsu73}).


According to the definition given in \eqref{defiRamifiedWeightFunction}, for any $t,s>0$ we get
\begin{equation}\label{omegaMspower}
(\omega_M)^s(t)=\omega_M(t^s)=\sup_{p\in\NN}\log\left(\frac{t^{sp}}{M_p}\right)= \sup_{p\in\NN}\log\left(\left(\frac{t^p}{(M_p)^{1/s}}\right)^s\right)=s\omega_{M^{1/s}}(t),
\end{equation}
where $M^{1/s}:=((M_p)^{1/s})_{p\in\NN}$.

We summarize some more well-known facts for this function:

\begin{lemma}\label{assofuncproper}
Let $M\in\hyperlink{LCset}{\mathcal{LC}}$.

\begin{itemize}

\item[$(i)$] $\omega_M$ is a normalized weight function satisfying \hyperlink{om3}{$(\omega_3)$} and \hyperlink{om4}{$(\omega_4)$}.

\item[$(ii)$]
$\lim_{p\rightarrow\infty}(m_p)^{1/p}=+\infty$ implies \hyperlink{om5}{$(\omega_5)$} for $\omega_M$.

\item[$(iii)$] $M$ has \hyperlink{mg}{$(\on{mg})$} if and only if $\omega_M$ has \hyperlink{om6}{$(\omega_6)$}.

\item[$(iv)$]
%
If $M$ satisfies
$\hyperlink{gamma1}{(\gamma_1)}$, then $\omega_M$ has \hyperlink{omsnq}{$(\omega_{\on{snq}})$}.


\end{itemize}
So for any strongly regular weight sequence $M$ the weight function $\omega_M$ satisfies \hyperlink{om3}{$(\omega_3)$}, \hyperlink{om4}{$(\omega_4)$} and \hyperlink{omsnq}{$(\omega_{\on{snq}})$}.

\end{lemma}

\demo{Proof}
$(i)$ See \cite[Definition 3.1]{Komatsu73}.

$(ii)$ That $\lim(m_p)^{1/p}=+\infty$ implies \hyperlink{om5}{$(\omega_5)$} for $\omega_M$ follows along the same lines as a similar argument in
\cite[Lemma 12 $(iv)\Rightarrow(v)$]{BonetMeiseMelikhov07}.

$(iii)$ See \cite[Proposition 3.6]{Komatsu73}.

$(iv)$
%
It follows from \cite[Proposition 4.4]{Komatsu73}.

\qed\enddemo

\begin{lemma}\label{assofuncproperbis}
Let $\omega$ be a normalized weight satisfying  \hyperlink{om3}{$(\omega_3)$}, and $\Omega=\{W^x=(W^x_p)_{p\in\NN}: x>0\}$ be its associated weight matrix, where $W^x_p=\exp\left(\frac{1}{x}\varphi^{*}_{\omega}(xp)\right)$.
Then, we get
    \begin{equation*}
    \forall\;x>0\;\forall\;t\ge 0:\;\;\;x\omega_{W^x}(t)\le\omega(t).
    \end{equation*}
    If $\omega$ satisfies moreover \hyperlink{om4}{$(\omega_4)$}, then $\omega\hyperlink{sim}{\sim}\omega_{W^x}$ for each $x>0$, more precisely we get
    \begin{equation}\label{goodequivalence}
    \forall\;x>0\;\exists\;C_x>0\;\forall\;t\ge 0:\;\;\;x\omega_{W^x}(t)\le\omega(t)\le 2x\omega_{W^x}(t)+C_x.
    \end{equation}
\end{lemma}

\demo{Proof}
For these estimates we recall \cite[Lemma 5.7]{compositionpaper}, respectively \cite[Theorem 4.0.3, Lemma 5.1.3]{dissertation}: Given the normalized weight $\omega$ satisfying \hyperlink{om3}{$(\omega_3)$} and the sequence $W^1=(\exp\left(\varphi^{*}_{\omega}(p)\right))_{p\in\NN}$, there exists some $c>0$ such that for all $t\ge 0$ we get $\omega_{W^1}(t)\le\omega(t)\le 2\omega_{W^1}(t)+c$, where the first inequality does not need \hyperlink{om4}{$(\omega_4)$} while the second does. Now, for any $x>0$ we will apply the previous statement to the weight $\tau_x(t):=\omega(t)/x$, which has the same properties assumed for $\omega$. We need to compute, for $p\in\NN$, the value
\begin{align*}
\varphi^{*}_{\tau_x}(p)&=\sup_{y\ge 0}\left\{py-\varphi_{\tau_x}(y)\right\}=\sup_{y\ge 0}\left\{py-\tau_x(e^y)\right\}=\sup_{y\ge 0}\left\{py-\frac{1}{x}\omega(e^y)\right\}\\
&=
\sup_{y\ge 0}\left\{py-\frac{1}{x}\varphi_{\omega}(y)\right\}=\frac{1}{x}\sup_{y\ge 0}\{(xp)y-\varphi_{\omega}(y)\}=\frac{1}{x}\varphi^{*}_{\omega}(xp).
\end{align*}
So, it turns out that, in the same way that $W^1$ was the sequence corresponding to $\omega$ in our previous statement, the sequence $W^x$ is the one corresponding to $\tau_x$, and so there exists some $c_x>0$ such that for all $t\ge 0$ we get $\omega_{W^x}(t)\le\tau_x(t)\le 2\omega_{W^x}(t)+c_x$, under the same conditions as before. The conclusion is immediate by the definition of $\tau_x$.

\qed\enddemo


Another important function will be introduced now: Let $M\in\RR_{>0}^{\NN}$ ($M_0=1$) and put
\begin{equation*}
h_M(t):=\inf_{k\in\NN}M_k t^k.
\end{equation*}
The functions $h_M$ and $\omega_M$ are related by
\begin{equation}\label{functionhequ2}
h_M(t)=\exp(-\omega_M(1/t)),\quad\;t>0,
\end{equation}
since $\log(h_M(t))=\inf_{k\in\NN}\log(t^kM_k)=-\sup_{k\in\NN}-\log(t^kM_k)=-\omega_M(1/t)$ (e.g. see also \cite[p. 11]{ChaumatChollet94}). By definition we immediately get:

\begin{lemma}\label{functionhproperties}
Let $M,N\in\RR_{>0}^{\NN}$ be given, then
\begin{itemize}
\item[$(i)$] The function $h_M(t)$ is nondecreasing,
\item[$(ii)$] If $M$ is normalized, then $h_{M}(t)\le 1$ for all $t>0$; if moreover $M\in\hyperlink{LCset}{\mathcal{LC}}$, then $h_M(t)=1$ for all $t$ sufficiently large and $\lim_{t\rightarrow 0}h_M(t)=0$,
\item[$(iii)$] $M\le N$ implies $h_M\le h_N$, more generally $M\hyperlink{precsim}{\precsim}N$ implies that $h_M(t)\le h_N(Ct)$ holds for some $C\ge 1$ and all $t>0$,
\end{itemize}
\end{lemma}

From Lemma~\ref{assofuncproperbis} and the equality~\eqref{functionhequ2} we immediately deduce the following result. Here, we write as before $\omega^{\iota}(t)=\omega(1/t)$ for a given weight function $\omega$.

\begin{lemma}\label{lemmaexponentialdecrease.hm}
Let $\omega$ be a normalized weight function satisfying \hyperlink{om3}{$(\omega_3)$}, then we get
\begin{equation}\label{equaExp-Omega.less.hm}
    \forall\;x>0\;\forall\;t\ge 0:\;\;\;\exp(-\omega^{\iota}(t))\le (h_{W^x}(t))^x.
    \end{equation}
If moreover $\omega$ has \hyperlink{om4}{$(\omega_4)$}, then
\begin{equation}\label{equaExp-Omega.bigger.hm}
    \forall\;x>0\;\exists\;C_x>0\;\forall\;t\ge 0:\;\;\;\exp(-C_x)(h_{W^x}(t))^{2x}\le \exp(-\omega^{\iota}(t)).
    \end{equation}
\end{lemma}

\section{Legendre conjugates and equivalent weight functions}\label{sectLegendreConjWeightMatrtoWeightFunc}

\subsection{Legendre conjugates of a weight $\omega$}\label{conugateofaweight}

In the study of ultradifferentiable or ultraholomorphic classes defined by weight sequences, the operations of multiplying or dividing the sequence by the factorials play a prominent role, and this paper will not be an exception. When weight functions are considered instead, the corresponding analogous operations are expressed by means of Legendre conjugates, which will be described now. Moreover, in Section~\ref{growthindices} suitable growth indices will be introduced for both weight sequences and functions, and it will be seen how these values are increased or decreased by 1 when applying the corresponding operations in both cases.

For any $M\in\hyperlink{LCset}{\mathcal{LC}}$ with $\lim_{p\rightarrow\infty}(m_p)^{1/p}=+\infty$ there exists a
close connection between $\omega_m$ and a conjugate for $\omega_M$,
as considered in \cite[Definition 1.4]{PetzscheVogt} and \cite{BonetBraunMeiseTaylorWhitneyextension}, see Lemma \ref{omegaconjugateequivalent} below for more details. This conjugate must not be mixed with $\varphi^{*}_{\omega}$ as considered in \eqref{legendreconjugate}. We warn the reader that the terminology differs from one author to another: In \cite{Beurling72} the conjugates are named after Legendre, in \cite{PetzscheVogt} after Young.

Let $\omega$
be a weight function, then for any $s\ge 0$ we define
\begin{equation*}
\omega^{\star}(s):=\sup_{t\ge 0}\{\omega(t)-st\}.
\end{equation*}
$\omega^{\star}$ is the \textit{upper Legendre conjugate (or upper Legendre envelope)} of $\omega$.

We summarize some basic properties, see also \cite[Remark 1.5]{PetzscheVogt}.
\begin{itemize}
\item[$(i)$] By definition, $\omega^{\star}(0)=+\infty$. If $\omega$ has in addition \hyperlink{om5}{$(\omega_5)$}, then $\omega^{\star}(s)<+\infty$ for all $s>0$: Indeed, we have that for any $s>0$ (however small) there exists some $C_s>0$ (large enough) such that for all $t\ge0$ we get $\omega(t)\le st+C_s$, and so $\omega^{\star}(s)\le C_s$. In this case, the function $\omega^{\star}:(0,+\infty)\rightarrow[0,+\infty)$ is nonincreasing, continuous and convex, and $\lim_{s\to 0}\omega^{\star}(s)=+\infty$, $\lim_{s\to \infty}\omega^{\star}(s)=0$.


\item[$(ii)$] So, whenever the weight function $\omega$ has \hyperlink{om5}{$(\omega_5)$},
the function $(\omega^{\star})^{\iota}$
is again a (generally not normalized) weight function if we set $(\omega^{\star})^{\iota}(0):=0$ (recall that $(\omega^{\star})^{\iota}(t)=\omega^{\star}(1/t)$, $t>0$).
\end{itemize}

We introduce now a new conjugate. For any $h:(0,+\infty)\rightarrow[0,+\infty)$ which is nonincreasing and such that $\lim_{s\rightarrow 0}h(s)=+\infty$, we can define the so-called {\itshape lower Legendre conjugate (or envelope)} $h_{\star}:[0,+\infty)\rightarrow[0,+\infty)$ of $h$ by
\begin{equation*}
h_{\star}(t):=\inf_{s>0}\{h(s)+ts\},\hspace{15pt}t\ge 0.
\end{equation*}
$h_{\star}$ is clearly nondecreasing, continuous and concave, and $\lim_{t\rightarrow\infty}h_{\star}(t)=\infty$, see \cite[(8), p. 156]{Beurling72}. Moreover, if $\lim_{s\to \infty}h(s)=0$ then $h_{\star}(0)=0$, and so $h_{\star}$ is a weight function.

In our work this second conjugate will be mainly applied to the case
$h(t):=\omega^{\iota}(t)=\omega(1/t)$, where $\omega$ is a weight function, so that $(\omega^{\iota})_{\star}$ is again a weight function; in particular, we will frequently find the case $h(t)=\omega^{\iota}_M(t)=\omega_M(1/t)$ for $M\in\RR_{>0}^{\NN}$ with $\lim_{p\rightarrow\infty}(M_p)^{1/p}=+\infty$.

In \cite[Proposition 1.6]{PetzscheVogt} it was shown that for any $\omega:[0,+\infty)\rightarrow[0,+\infty)$ concave and nondecreasing we get
\begin{equation*}
\forall\;t>0:\;\;\;\;\omega(t)=\inf_{s>0}\{\omega^{\star}(s)+st\}= (\omega^{\star})_{\star}(t).
\end{equation*}
In case $\omega$ is a weight function satisfying \hyperlink{om5}{$(\omega_5)$}, $(\omega^{\star})_{\star}$ is a weight function and it is indeed the least concave majorant of $\omega$ (in the sense that, if $\tau:[0,+\infty)\rightarrow[0,+\infty)$ is concave and $\omega\le\tau$, then $(\omega^{\star})_{\star}\le\tau$), see~\cite{roever}.

We prove now several properties for $\omega^{\star}$ which will be needed below. For this we use $0^0:=1$ and recall the following consequence of {\itshape Stirling's formula:}
\begin{equation}\label{stirling}
\forall\;n\in\NN:\;\;\;\left(\frac{n}{e}\right)^n\le n!\le n^n.
\end{equation}

\begin{lemma}\label{omegaconjugateequivalent}
\begin{itemize}
\item[$(i)$] Let $\sigma$ and $\tau$ be two weight functions with \hyperlink{om5}{$(\omega_5)$}, and suppose there exist $A,B>0$ such that $$\forall\;t>0:\;\;\;\tau(t)\le A\sigma(t)+B.$$
    Then
$$\forall\;s>0:\;\;\;\tau^{\star}(s)\le A\sigma^{\star}\left(\frac{s}{A}\right)+B.$$
Consequently $\sigma\hyperlink{sim}{\sim}\tau$ implies
\begin{equation*}
\exists\;C\ge 1\;\forall\;s>0:\;\;\;-C+C^{-1}\sigma^{\star}(Cs)\le\tau^{\star}(s)\le C\sigma^{\star}\left(\frac{s}{C}\right)+C.
\end{equation*}

\item[$(ii)$] Let $M\in\RR_{>0}^{\NN}$ such that $\lim_{p\rightarrow\infty}(p!^{1-b}m_p)^{1/p}=+\infty$ for some $b>0$. Then for each $0<a\le b$ the mappings $s\mapsto((\omega_{M})^a)^{\star}(s)$ and $s\mapsto\omega_{M/G^a}(1/s^a)$, with $G^a_p:=p!^a$, are equivalent in the previous sense, more precisely
\begin{equation}\label{Dynkinequiv0}
\forall\;s>0:\;\;\;((\omega_{M})^a)^{\star}(s)\le\omega_{M/G^a}\left(\frac{a^a}{s^a}\right)\le((\omega_{M})^a)^{\star}\left(\frac{s}{e}\right).
\end{equation}

\end{itemize}
\end{lemma}

\demo{Proof}
$(i)$ 
Let $s>0$, then
$$\tau^{\star}(s)=\sup_{t\ge 0}\{\tau(t)-st\}\le\sup_{t\ge 0}\{A\sigma(t)+B-st\}=A\sup_{t\ge 0}\{\sigma(t)-(sA^{-1})t\}+B=A\sigma^{\star}(\frac{s}{A})+B,$$
see also \cite[Remark 1.7]{PetzscheVogt}.

$(ii)$ For any $a\in(0,b]$ we have $\lim_{p\rightarrow\infty}(M_p/(p!)^a)^{1/p}= \lim_{p\rightarrow\infty}(p!^{1-a}m_p)^{1/p}=+\infty$ if, and only if, $\lim_{p\rightarrow\infty}(M^{1/a}_p/p!)^{1/p}= \lim_{p\rightarrow\infty}(p!^{1-a}m_p)^{1/(ap)}=+\infty$,
hence we may apply Lemma~\ref{assofuncproper}(ii) to deduce that $\omega_{M^{1/a}}$ satisfies \hyperlink{om5}{$(\omega_5)$}. 
Consequently, as indicated in the study of the properties of the upper Legendre conjugate, the function $(\omega_{M^{1/a}})^{\star}$ is well-defined from $(0,\infty)$ to $(0,\infty)$, and by \eqref{omegaMspower} coincides with $(a(\omega_{M})^a)^{\star}$.

We follow now the proof \cite[Lemma 5.7.8]{diplomadebrouwere}, where only the case $a=1$ was treated. Let $s>0$, then
\begin{align*}
(a^{-1}\omega_{M^{1/a}})^{\star}(s)&=(\omega^a_{M})^{\star}(s):= \sup_{t\ge 0}\{\omega_M(t^a)-st\}=\sup_{t\ge 0}\left\{\sup_{p\in\NN}\log\left(\frac{t^{ap}}{M_p}\right)-st\right\}\\
&=\sup_{p\in\NN}\sup_{t\ge 0}\left\{\log\left(\frac{t^{ap}}{M_p}\right)-st\right\}.
\end{align*}
For $s>0$ and $p\in\NN$ fixed we consider $f_{s,p}:(0,+\infty)\rightarrow\RR$ defined by $$f_{s,p}(t):=\log\left(\frac{t^{ap}}{M_p}\right)-st= ap\log(t)-\log(M_p)-st,\quad p\in\NN_{>0}; \ f_{s,0}(t)=-st.
$$
Hence, $\sup_{t> 0}f_{s,0}=0$ for any $s>0$. Let $p\ge 1$, then $f'_{s,p}(t)=ap\frac{1}{t}-s=0\Leftrightarrow t=\frac{ap}{s}$ and $f_{s,p}$ attains its maximum at this point. Hence we get $$f_{s,p}(ap/s)=ap\log(ap/s)-\frac{ap}{s}s-\log(M_p)=\log\left(\frac{a^{ap}p^{ap}}{s^{ap}M_p}\right)-\log(\exp(ap))=\log\left(\left(\frac{a^a}{s^a}\right)^p\frac{p^{ap}}{e^{ap}M_p}\right),$$
which holds also for $p=0$ by $0^0:=1$. Thus we have shown
$$\forall\;s>0:\;\;\;((\omega_{M})^a)^{\star}(s)=\sup_{p\in\NN}\log\left(\left(\frac{a^a}{s^a}\right)^p\frac{p^{ap}}{e^{ap}M_p}\right).$$
The left hand side of \eqref{stirling} gives $\frac{p^{ap}}{(e^a)^pM_p}\le\frac{p!^a}{M_p}$ for all $p\in\NN$, i.e. the left hand side of \eqref{Dynkinequiv0}. By the right hand side of \eqref{stirling} we get $\log\left(\frac{p!^a}{M_p}\right)\le\log\left(\frac{p^{ap}}{M_p}\right)$ for all $p\in\NN$ and so the right hand side of \eqref{Dynkinequiv0}.

\qed\enddemo

Combining the previous lemma with results from \cite[Section 5]{compositionpaper} we get the following consequences, which have already appeared, in a weaker form, in~\cite{whitneyextensionmixedweightfunction}.

\begin{corollary}\label{omegaconjugateequivalentcor}
Let $\omega$ be a normalized weight with \hyperlink{om3}{$(\omega_3)$}, \hyperlink{om4}{$(\omega_4)$} and \hyperlink{om5}{$(\omega_5)$}, let $\Omega=\{W^x=(W^x_p)_{p\in\NN}: x>0\}$ be its associated weight matrix, and put $w^x=(W^x_p/p!)_{p\in\NN}$, $x>0$. Then,
\begin{equation}\label{omegaconjugateequivalentcorequ1}
\forall\;x>0\;\exists\;C_x\ge 1\;\forall\;s>0:\;\;\;x\omega_{W^x}^{\star}(\frac{s}{x})\le\omega^{\star}(s)\le 2x\omega_{W^x}^{\star}(\frac{s}{2x})+C_x
\end{equation}
and
\begin{equation}\label{omegaconjugateequivalentcorequ2}
\forall\;x>0\;\exists\;C_x\ge 1\;\forall\;s>0:\;\;\;x\omega_{w^x}\left(\frac{x}{es}\right)\le\omega^{\star}(s)\le 2x\omega_{w^x}\left(\frac{2x}{s}\right)+C_x,
\end{equation}
or equivalently,
\begin{equation*}
\forall\;x>0\;\exists\;C_x\ge 1\;\forall\;s>0:\;\;\;h_{w^x}\left(\frac{es}{x}\right)^x\ge\exp(-\omega^{\star}(s))\ge \exp(-C_x)h_{w^x}\left(\frac{s}{2x}\right)^{2x},
\end{equation*}
where, for all the inequalities on the left to hold, it is not necessary to impose \hyperlink{om4}{$(\omega_4)$}.
\end{corollary}

\demo{Proof}
To prove \eqref{omegaconjugateequivalentcorequ1} we apply \eqref{goodequivalence}, the stability of \hyperlink{om5}{$(\omega_5)$} under equivalence, and Lemma \ref{omegaconjugateequivalent}(ii). For \eqref{omegaconjugateequivalentcorequ2} we depart from \eqref{omegaconjugateequivalentcorequ1} and recall (see Subsection~\ref{weightmatrixfromfunction}) that \hyperlink{om5}{$(\omega_5)$} implies   $\lim_{p\rightarrow\infty}(w^x_p)^{1/p}=+\infty$ for every $x>0$, so we may apply also $(ii)$ in the previous result for $a=1$. The last inequalities are just a re-writing of \eqref{omegaconjugateequivalentcorequ2} thanks to the very definition \eqref{functionhequ2}.
\qed\enddemo

\subsection{Equivalent weight functions from Legendre conjugates}

The results in this subsection state the equivalence of different weight functions, naturally appearing in our arguments and in whose construction the Legendre transformations play a significant role.

\begin{remark}\label{remarkequivweight}
Let $M,N\in\hyperlink{LCset}{\mathcal{LC}}$ 
and let $\omega_M$ have \hyperlink{om1}{$(\omega_1)$}. Then $M\hyperlink{approx}{\approx}N$ implies $\omega_M\hyperlink{sim}{\sim}\omega_N$ as follows: First $M\hyperlink{approx}{\approx}N$ implies $\omega_M(A^{-1}t)\le\omega_N(t)\le\omega_M(At)$ for all $t\ge 0$ and some $A\ge 1$, hence by iterating \hyperlink{om1}{$(\omega_1)$} there exists some $B\ge 1$ such that for all $t$ sufficiently large:
$$B^{-1}\omega_M(t)\le\omega_M(A^{-1}t)\le\omega_N(t)\le\omega_M(At)\le B\omega_M(t).$$
\end{remark}

Moreover we recall (see Lemma~\ref{omegaconjugateequivalent}(ii) with $a=1$) that for $M\in\RR_{>0}^{\NN}$ with $\lim_{p\rightarrow\infty}(m_p)^{1/p}=+\infty$ one has
\begin{equation}\label{Dynkinequiv1}
\forall\;s>0:\;\;\;\omega_{M}^{\star}(s)\le
\omega_{m}\left(\frac{1}{s}\right)\le
\omega_{M}^{\star}\left(\frac{s}{e}\right).
\end{equation}

\begin{lemma}\label{omegaconjugateequivalentbis}
\begin{itemize}
\item[$(i)$] Let $M\in\RR_{>0}^{\NN}$ with $\lim_{p\rightarrow\infty}(m_p)^{1/p}=+\infty$ be given and define the sequence
\begin{equation*}
N_p:=\sup_{t>0}\frac{t^p}{\exp((\omega^{\iota}_m)_{\star}(t))},
\end{equation*}
where $\omega^{\iota}_m(t)=\omega_m(1/t)$. Then we get $N\hyperlink{approx}{\approx}(p!m^{\on{lc}}_p)_p$ (see~\eqref{equaLCminorant} for the definition of $(m^{\on{lc}}_p)_p$),
and so $N$ is log-convex and equivalent to a strongly log-convex sequence.

\item[$(ii)$] Let $Q\in\RR_{>0}^{\NN}$ such that $Q\hyperlink{approx}{\approx}M$ and $m\in\hyperlink{LCset}{\mathcal{LC}}$. Then $\omega_Q$ is equivalent to a concave function, more precisely we get
\begin{equation}\label{Dynkinequiv3}
\forall\;x\ge\mu_1:\;\;\;\omega_M(x)\le(\omega^{\iota}_m)_{\star}(x)\le 1+\omega_M(ex),
\end{equation}
and since $(\omega^{\iota}_m)_{\star}$ is concave, we have $\omega_Q\hyperlink{sim}{\sim}(\omega^{\iota}_m)_{\star}$.
\end{itemize}
\end{lemma}

\demo{Proof}
$(i)$ For all $p\in\NN$ we get
\begin{align*}
\sup_{t>0}\frac{t^p}{\exp((\omega^{\iota}_m)_{\star}(t))}&= \exp\left(\sup_{t>0}\{p\log(t)-(\omega^{\iota}_m)_{\star}(t)\}\right)\\
&=\exp\left(\sup_{t>0}\{p\log(t)-\inf_{s>0}\{\omega_m(1/s)+st\}\}\right)\\
&=\exp\left(\sup_{t,s>0}\{p\log(t)-\omega_m(1/s)-st\}\right).
\end{align*}
Let $p\in\NN$ and $s>0$ be fixed and put $$f_{p,s}(t):=p\log(t)-\omega_m(1/s)-st,\ \ p\in\NN_{>0};\quad f_{0,s}(t)=-\omega_m(1/s)-st.
$$
Clearly, $\sup_{t>0}f_{0,s}=-\omega_m(1/s)$. For all $p\ge 1$ we get $f'_{p,s}(t)=\frac{p}{t}-s=0\Leftrightarrow t=\frac{p}{s}$, the point where $f_{p,s}$ attains its maximum. Hence $$f_{p,s}(p/s)=p\log(p/s)-\omega_m(1/s)-p= \log\left(\frac{p^p}{(es)^p}\right)-\omega_m(1/s),
$$
which holds also for the case $p=0$ by $0^0:=1$. Thus for any $p\in\NN$,
\begin{align*}
\sup_{t>0}\frac{t^p}{\exp((\omega^{\iota}_m)_{\star}(t))}&= \exp\left(\sup_{s>0}\{\log\left(\frac{p^p}{(es)^p}\right)- \omega_m(1/s)\}\right)\\
&=\frac{p^p}{e^p}\sup_{s>0}\frac{1}{s^p\exp(\omega_m(1/s))}
=\frac{p^p}{e^p}\sup_{s>0}\frac{s^p}{\exp(\omega_m(s))}= \frac{p^p}{e^p}m^{\on{lc}}_p,
\end{align*}
where in the last step we have applied~\eqref{equaLCminorant}.\vspace{6pt}

$(ii)$ We follow and recall the arguments of \cite[p. 233]{normality}. On the one hand we get by \eqref{Dynkinequiv1} for any $x\ge 0$:
\begin{align*}
(\omega_m^{\iota})_{\star}(x)&=\inf_{y>0}\{\omega_m(1/y)+xy\}\ge\inf_{y>0}\{\omega^{\star}_M(y)+xy\}=\inf_{y>0}\{\sup_{u\ge 0}\{\omega_M(u)-uy\}+xy\}
\\&
=\inf_{y>0}\sup_{u\ge 0}\{\omega_M(u)+y(x-u)\}\underbrace{\ge}_{x=u}\omega_M(x).
\end{align*}
For this estimate the log-convexity of $m$ was not used. On the other hand, first we get
\begin{align*}
\exp(-(\omega_m^{\iota})_{\star}(x))&=\exp\left(-\inf_{y>0}\{\omega_m(1/y)+xy\}\right)=\exp\left(\sup_{y>0}\{-xy-\omega_m(1/y)\}\right)
\\&
=\sup_{y>0}\{\exp(-xy)\exp(-\omega_m(1/y))\}=\sup_{y>0}\{\exp(-xy)h_m(y)\}.
\end{align*}
For convenience, we write $\mu^{*}_{0}:=1$, $\mu^{*}_{n}:=\mu_n/n$, $n\in\NN_{>0}$. Let now $x\in[\mu_n,\mu_{n+1})$, $n\ge 1$, and put $y_0:=\max\{\frac{n}{x},\frac{1}{\mu^{*}_{n+1}}\}$. Hence $n\frac{1}{\mu_{n+1}}<\frac{n}{x}\le n\frac{1}{\mu_n}$ and by the strong log-convexity $\frac{1}{\mu^{*}_{n+1}}\le n\frac{1}{\mu_n}=\frac{1}{\mu^{*}_n}$. So $\frac{1}{\mu^{*}_{n+1}}\le y_0\le\frac{1}{\mu^{*}_n}$ and we estimate as follows:
$$\exp(-(\omega_m^{\iota})_{\star}(x))= \sup_{y>0}\{\exp(-xy)h_m(y)\}\ge\exp(-xy_0)h_m(y_0)=m_ny_0^n\exp(-xy_0),$$
where the last equality holds by the choice of $y_0$ as explained above. Finally, we have to show that $e\exp(\omega_M(ex))=e\sup_{l\in\NN}\frac{e^lx^l}{M_l}\ge\frac{1}{m_ny_0^n}\exp(xy_0)$. It suffices to consider the choice $l=n$ on the left hand side which yields $\frac{e^{n+1}x^n}{n!}\ge\frac{1}{y_0^n}\exp(xy_0)$. And this holds true since, on the one hand, $\frac{1}{y_0^n}\le\frac{x^n}{n^n}\le\frac{x^n}{n!}$, and on the other hand, $\exp(x(n/x))=\exp(n)$ and $\exp(x/\mu_{n+1}^{*})\le\exp(\mu_{n+1}/\mu^{*}_{n+1})=\exp(n+1)$,  which together proves $\exp(xy_0)\le\exp(n+1)$.\vspace{6pt}

So far we have shown \eqref{Dynkinequiv3}.  $(\omega_m^{\iota})_{\star}$ has \hyperlink{om1}{$(\omega_1)$} by concavity, and so we deduce that $\omega_M\hyperlink{sim}{\sim}(\omega_m^{\iota})_{\star}$ and that also $\omega_M$ has \hyperlink{om1}{$(\omega_1)$}. Finally, since $Q\hyperlink{approx}{\approx}M$, Remark \ref{remarkequivweight} yields $\omega_Q\hyperlink{sim}{\sim}\omega_M$ and we are done.
\qed\enddemo

Using this result we can prove the following Corollary:

\begin{corollary}\label{omegaconjugateequivalentcor1}
Let $M\in\RR_{>0}^{\NN}$ with $\lim_{p\rightarrow\infty}(m_p)^{1/p}=+\infty$. Then the functions $\omega_L$, where $L_p:=p!m^{\on{lc}}_p$, and $(\omega_m^{\iota})_{\star}$ are equivalent with respect to $\hyperlink{sim}{\sim}$.
\end{corollary}

\demo{Proof}
Recall that, by the very definition of $(\omega_m^{\iota})_{\star}$, we have $(\omega_m^{\iota})_{\star}=(\omega_{m^{\on{lc}}}^{\iota})_{\star}$, and so, by $(ii)$ in Lemma \ref{omegaconjugateequivalentbis} applied for $Q=L$, i.e. with $m^{\on{lc}}$ instead of $m$, we get $(\omega_m^{\iota})_{\star}\hyperlink{sim}{\sim}\omega_L$.
\qed\enddemo

The following result will be useful later on.
\begin{lemma}\label{lemmaWeightequivConcave}
Let $\omega$ be a normalized weight function satisfying \hyperlink{om3}{$(\omega_3)$}, \hyperlink{om4}{$(\omega_4)$} and \hyperlink{om5}{$(\omega_5)$}, and which is moreover equivalent to a concave weight function. Let $\Omega=\{W^x:x>0\}$ be the weight matrix associated with $\omega$, and
define $L_p^x:=p!(w_p^x)^{\on{lc}}$. Then, for all $x>0$ it holds $\omega_{W^x}\hyperlink{sim}{\sim}\omega_{L^x}$.
\end{lemma}
\demo{Proof}
As indicated in Subsection~\ref{weightmatrixfromfunction}.(iv),  \hyperlink{om5}{$(\omega_5)$} for $\omega$ implies that $\lim_{p\to\infty}(w_p^x)^{1/p}=\infty$, so the function $\omega_{w^x}$ is well-defined.
Note that, as $\omega\hyperlink{sim}{\sim}\omega_{W^x}$, $\omega_{W^x}$ is equivalent to a concave weight function, and it will also be equivalent to its least concave majorant, $\kappa^x:=(\omega_{W^x}^{\star})_{\star}$.
Taking into account \eqref{Dynkinequiv1}, for every $t>0$ we deduce that
\begin{align}\label{equaEquivSigma_xOmega_wx}
\kappa^x(t)&=\inf_{s>0}\{\omega_{W^x}^{\star}(s)+ts\}\le \inf_{s>0}\{\omega_{w^x}^{\iota}(s)+ts\}\nonumber\\
&=(\omega_{w^x}^{\iota})_{\star}(t)\le \inf_{s>0}\{\omega_{W^x}^{\star}(s/e)+ts\}=\inf_{s>0}\{\omega_{W^x}^{\star}(s)+ets\}
=\kappa^x(et).
\end{align}
The concavity of $\kappa^x$ implies the existence of a constant $C>0$ such that $\kappa^x(et)\le C\kappa^x(t)$ for every $t>0$, and so the inequalities~\eqref{equaEquivSigma_xOmega_wx} guarantee that $(\omega_{w^x}^{\iota})_{\star}\hyperlink{sim}{\sim}\kappa^x$. Finally,
$$
\omega_{W^x}\hyperlink{sim}{\sim}\kappa^x\hyperlink{sim}{\sim}(\omega_{w^x}^{\iota})_{\star}\hyperlink{sim}{\sim}\omega_{L^x},
$$
where the last equivalence is due to Corollary~\ref{omegaconjugateequivalentcor1}.
\qed\enddemo

\section{Growth indices}\label{growthindices}
\subsection{The growth index $\gamma(M)$ introduced by V. Thilliez}\label{thilliezgrowthindex}

We revisit the definition of the growth index $\gamma(M)$ introduced in \cite[Section 1.3]{Thilliezdivision}. This is necessary if we pretend to explain the result about the mixed setting as a complement to the extension results by V.~Thilliez.

Let $\gamma\in\RR$ be given, then $M\in\RR_{>0}^{\NN}$ satisfies $(P_{\gamma})$ (see \cite[Definition 1.3.1]{Thilliezdivision} where only $\gamma>0$ was considered), if
\begin{equation}\label{Pgamma1}
\exists\;\nu=(\nu_p)_p\;\exists\;a\ge 1\;\forall\;p\in\NN:\;\;\;a^{-1}\mu_p\le\nu_p\le a\mu_p
\end{equation}
and such that
\begin{equation}\label{Pgamma2}
p\mapsto\frac{\nu_p}{p^{\gamma}}\;\;\text{is nondecreasing}.
\end{equation}
\eqref{Pgamma1} is precisely $M\hyperlink{simeq}{\simeq}N$ for $N=(N_p)_p$ given by $N_0:=1$, $N_p=\prod_{j=1}^p \nu_j$, $p\in\NN_{>0}$.
\eqref{Pgamma2} is equivalent to the fact that $(\frac{N_p}{(p!)^{\gamma}})_{p\in\NN}$
is log-convex.  The growth index of $M$, introduced in \cite[Definition 1.3.5]{Thilliezdivision}, is defined by
\begin{equation*}
\gamma(M):=\sup\{\gamma\in\RR: (P_{\gamma})\;\;\text{is satisfied}\}.
\end{equation*}
If $\{\gamma\in\RR: (P_{\gamma})\;\;\text{is satisfied}\}=\emptyset$, then we put $\gamma(M):=-\infty$, if $\{\gamma\in\RR: (P_{\gamma})\;\;\text{is satisfied}\}=\RR$, then $\gamma(M):=+\infty$. We point out that Thilliez only considered the case $M\in\hyperlink{LCset}{\mathcal{LC}}$, and so $\gamma(M)\ge 0$.
We summarize some properties for $\gamma(M)$:

Note that $m$ has $(P_{\gamma})$ if, and only if, $M$ has $(P_{\gamma+1})$ (recall that $M_p=p!m_p$ for all $p\in\NN$). Then, by definition $\gamma(m)+1=\gamma(M)$ holds.


Moreover, $M\hyperlink{simeq}{\simeq}N$ implies $\gamma(M)=\gamma(N)$.



    In \cite[Lemma 1.3.2]{Thilliezdivision} it was shown that for $m\in\hyperlink{LCset}{\mathcal{LC}}$ such that $M$ satisfies \hyperlink{gamma1}{$(\gamma_1)$} we always have $\gamma(m)>0$. For this implication the assumption $M\in\hyperlink{LCset}{\mathcal{LC}}$ is sufficient since we use \cite[Corollary 1.3]{petzsche} and which implies that always $(P_{\gamma})$ is satisfied for $m$ for some $\gamma>0$.



Combining several results~(\cite{BariSteckin}, \cite[Corollary 1.3]{petzsche}, \cite[Lemma 4.5]{Tikhonov04}) we may obtain the following useful information relating the condition \hyperlink{gamma1}{$(\gamma_1)$} to the value of the index $\gamma(M)$. A detailed proof will be included in~\cite{JimenezSanzSchindlIndices}.

\begin{lemma}\label{Thilliezstrongcondition}
Let $M\in\RR_{>0}^{\NN}$ be given, the following are equivalent::

\begin{itemize}
\item[$(i)$] $\gamma(M)>1$,

\item[$(ii)$] there exists $n\in\hyperlink{LCset}{\mathcal{LC}}$, $n\hyperlink{simeq}{\simeq}m$, and $N$ has \hyperlink{gamma1}{$(\gamma_1)$},

\item[$(iii)$] there exists $N\in\hyperlink{LCset}{\mathcal{LC}}$, $N\hyperlink{simeq}{\simeq}M$, and $N$ has \hyperlink{gamma1}{$(\gamma_1)$}.
\end{itemize}
In particular, if $M\hyperlink{simeq}{\simeq}N$ with $N\in\hyperlink{LCset}{\mathcal{LC}}$, this yields the following equivalence (which should be compared with Lemma \ref{strongweightspace} in the weight function setting):
\begin{itemize}
\item[$(+)$] $\gamma(M)>1$,

\item[$(++)$] $M$ satisfies \hyperlink{gamma1}{$(\gamma_1)$}.
\end{itemize}
\end{lemma}

\subsection{Growth index $\gamma(\omega)$}\label{growthindexgamma}

In this paragraph we introduce a growth index for a (not necessarily normalized) weight function. We are inspired by the equivalence $(ii)\Leftrightarrow(iii)$ in Proposition~\ref{Prop13MT88}.
It will turn out that some basic properties for a weight function are easily expressed in terms of this new index.

Let a weight function $\omega$ and $\gamma>0$ be given, we introduce the property
\begin{equation*}
(P_{\omega,\gamma}):\Longleftrightarrow\exists\;K>1:\;\;\;\limsup_{t\rightarrow\infty}\frac{\omega(K^{\gamma}t)}{\omega(t)}<K.
\end{equation*}
We note that if $(P_{\omega,\gamma})$ holds for some $K>1$, then also $(P_{\omega,\gamma'})$ is satisfied for all $\gamma'\le\gamma$ with the same $K$. Moreover we restrict ourselves to $\gamma>0$, because for $\gamma\le 0$ condition $(P_{\omega,\gamma})$ is satisfied for any weight $\omega$ (since it is nondecreasing and $K>1$).

Finally, we put
\begin{equation*}
\gamma(\omega):=\sup\{\gamma>0: (P_{\omega,\gamma})\;\;\text{is satisfied}\}.
\end{equation*}
So for any $0<s<\gamma(\omega)$ the weight $\omega^s$ given by $\omega^s(t)=\omega(t^s)$ has property \hyperlink{omsnq}{$(\omega_{\on{snq}})$}.

Let $\omega,\sigma$ satisfy $\sigma\hyperlink{sim}{\sim}\omega$ (or, equivalently, $\sigma^s\hyperlink{sim}{\sim}\omega^s$ for some $s>0$), then $\gamma(\sigma)=\gamma(\omega)$: Observe that each $(P_{\cdot,\gamma})$ is stable with respect to $\hyperlink{sim}{\sim}$ since \hyperlink{omsnq}{$(\omega_{\on{snq}})$} is clearly stable with respect to this relation. By definition and \eqref{omegaMspower} we immediately get
\begin{equation}\label{newindex3}
\forall\;s>0:\;\;\;\gamma(\omega^{1/s})=s\gamma(\omega).
\end{equation}

A first interesting result, whose proof will appear in~\cite{JimenezSanzSchindlIndices}, is the following.

\begin{lemma}\label{lemmagammaomegapositiveimpliesomega1}
Let $\omega$ be a weight function. Then, $\gamma(\omega)>0$ if, and only if, $\omega$ has \hyperlink{om1}{$(\omega_{1})$}.
\end{lemma}


Note that, while \hyperlink{om1}{$(\omega_{1})$} is a qualitative property of $\omega$, the condition $\gamma(\omega)>0$ is quantitative in the sense that the value of the index, as it will be shown in the next sections, provides an upper bound (except for the factor $\pi$) for the opening of the sectors in which extension results will be available for ultraholomorphic classes associated with $\omega$.

For a thorough study of the index $\gamma(\omega)$, its relationship with different properties for $\omega$ and the link with Thilliez's index $\gamma(M)$, we refer also to
~\cite{JimenezSanzSchindlIndices}. In this work we will need the following result, whose proof is included for completeness.

%
\begin{lemma}\label{strongweightspace}
$\omega$ satisfies \hyperlink{omsnq}{$(\omega_{\on{snq}})$} if and only if $\gamma(\omega)>1$.
\end{lemma}

\demo{Proof}
If $\gamma(\omega)>1$, then $\omega$ has \hyperlink{omsnq}{$(\omega_{\on{snq}})$} since $(P_{\omega,1})$ holds true (see Proposition~\ref{Prop13MT88}).

On the other hand let $\omega$ be given with \hyperlink{omsnq}{$(\omega_{\on{snq}})$}. Then, as already shown in \cite[Corollary 1.4]{MeiseTaylor88} there exists some $K>1$, $0<\alpha<1$ and $m\ge 0$ such that $\frac{\omega(Kt)}{\omega(t)}\le K^{\alpha}$ for all $t\ge K^m$. Take some $\beta$ such that $\alpha<\beta<1$ and so $\limsup_{t\rightarrow\infty}\frac{\omega(Kt)}{\omega(t)}\le K^{\alpha}<K^{\beta}$ is valid, what proves $(P_{\omega,\beta^{-1}})$ for any $\alpha<\beta<1$.
\qed\enddemo

Moreover we will have to use the following three results (for their proof we refer again to~\cite{JimenezSanzSchindlIndices}):

\begin{lemma}\label{growthindexgammalemma}
Let $\omega$ be a weight function. If $\gamma(\omega)>1$, then $(\omega^{\star})^{\iota}$
has \hyperlink{om1}{$(\omega_1)$} and we obtain
$$\gamma(\omega)=\gamma((\omega^{\star})^{\iota})+1.$$
In particular, $\gamma(\omega)>2$ implies that $(\omega^{\star})^{\iota}$ satisfies \hyperlink{omsnq}{$(\omega_{\on{snq}})$}.
\end{lemma}

\begin{corollary}\label{growthindexgammacorollary}
Let $M\in\hyperlink{LCset}{\mathcal{LC}}$ be given satisfying $\gamma(\omega_M)>1$. Then we obtain
$$\gamma(\omega_M)=\gamma(\omega_m)+1,$$
and consequently $\omega_m$ satisfies \hyperlink{om1}{$(\omega_1)$}.
In particular, $\gamma(\omega_M)>2$ implies that $\omega_m$ has \hyperlink{omsnq}{$(\omega_{\on{snq}})$}.

If $\omega$ is a normalized weight function with \hyperlink{om3}{$(\omega_3)$}, \hyperlink{om4}{$(\omega_4)$} and $\gamma(\omega)>1$, then $\gamma(\omega)=\gamma(\omega_{W^x})=\gamma(\omega_{w^x})+1=\gamma((\omega^{\star})^{\iota})+1$ for all $x>0$.
\end{corollary}

%
%



\begin{lemma}\label{lemmmaIndexLowerLegendreEnvelope}
For any weight function $\tau$ with $\gamma(\tau)>0$ one has
$$
\gamma((\tau^{\iota})_{\star})=\gamma(\tau)+1.
$$
\end{lemma}

%
%

\begin{remark}\label{remarkComparisonIndicesGamma}
We also mention without proof (see~\cite{JimenezSanzSchindlIndices}) that for a sequence $L\in\hyperlink{LCset}{\mathcal{LC}}$ one always has $\gamma(\omega_L)\ge \gamma(L)$, and that equality holds whenever $L$ has \hyperlink{mg}{(\text{mg})}. These facts will be useful for an extension result in a mixed setting that will be described in the last section of this paper.
\end{remark}


\section{Weight functions and matrices defining the same ultraholomorphic classes}

We will show now that, in two different situations, the ultraholomorphic class defined by a suitable weight function and the one associated with a related weight matrix coincide. These results will be important in order to prove the extension results we are looking for in a form convenient for the use of a truncated Laplace transform.

In the first situation, dealt with in the forthcoming three results, we assume that\par
\centerline{\textbf{$\tau$ is a normalized weight function
satisfying \hyperlink{om3}{$(\omega_3)$}, \hyperlink{om4}{$(\omega_4)$} and with $\gamma(\tau)>0$}.}
We recall that $\gamma(\tau)>0$ amounts to \hyperlink{om1}{$(\omega_1)$} (see Lemma~\ref{lemmagammaomegapositiveimpliesomega1}).

Denote by $\mathcal{T}:=\{T^x: x>0\}$ the associated weight matrix, i.e. $T^x_p:=\exp\left(\frac{1}{x}\varphi^{*}_{\tau}(xp)\right)$, and write also $\widehat{\mathcal{T}}:=\{\widehat{T}^x: x>0\}$, defined by $\widehat{T}^x_p:=p!T^x_p$ for each $x>0$ and $p\in\NN$.

\begin{lemma}\label{necessarytransformlemma1}
For all $x,y>0$ we get $(\omega^{\iota}_{T^x})_{\star}\hyperlink{sim}{\sim} (\omega^{\iota}_{T^y})_{\star}$.
\end{lemma}

\demo{Proof}
Let $x,y>0$ be arbitrary but fixed. By \cite[Lemma 5.7]{compositionpaper} we obtain \begin{equation*}
\omega_{T^x}\hyperlink{sim}{\sim}\tau\hyperlink{sim}{\sim}\omega_{T^y},
\end{equation*}
i.e. there exists some $C\ge 1$
such that $-C+C^{-1}\omega_{T^y}(s)\le\omega_{T^x}(s)\le C\omega_{T^y}(s)+C$ for all $s\ge 0$. Hence for any $s\ge 0$ we get:
\begin{align*}
(\omega^{\iota}_{T^x})_{\star}(s)&= \inf_{u>0}\left\{\omega_{T^x}\left(\frac{1}{u}\right)+us\right\}\le \inf_{u>0}\left\{C\omega_{T^y}\left(\frac{1}{u}\right)+us\right\}+C
\\&
=C\inf_{u>0}\left\{\omega_{T^y}\left(\frac{1}{u}\right)+ \frac{u}{C}s\right\}+C=C(\omega^{\iota}_{T^y})_{\star}(s/C)+C.
\end{align*}
Taking into account that each $(\omega^{\iota}_{T^x})_{\star}$ has \hyperlink{om1}{$(\omega_1)$} (by concavity), we have shown $(\omega^{\iota}_{T^x})_{\star}\hyperlink{sim}{\sim}(\omega^{\iota}_{T^y})_{\star}$ for all $x,y>0$.
\qed\enddemo

Combining previous results
we have shown so far:

\begin{corollary}\label{necessarytransformlemma1corollary}
The associated functions $\omega_{\widehat{T}^x}$ and $\omega_{\widehat{T}^y}$ are all equivalent with respect to $\hyperlink{sim}{\sim}$, more precisely $\omega_{\widehat{T}^x}\hyperlink{sim}{\sim}
(\omega_{T^y}^{\iota})_{\star}$ holds for all $x,y>0$.
\end{corollary}

\demo{Proof}
Let $x>0$, arbitrary but fixed, then take $m=T^x(=m^{\on{lc}})$ in Corollary \ref{omegaconjugateequivalentcor1} to show that $\omega_{\widehat{T}^x}\hyperlink{sim}{\sim}(\omega_{T^x}^{\iota})_{\star}$, what leads to the conclusion by Lemma \ref{necessarytransformlemma1}.
\qed\enddemo

Finally we can prove the following:

\begin{theorem}\label{theoEqualSpacesMatrixwithfactorialWeight}
For the considered weight $\tau$, the following identities hold as locally convex vector spaces for all sector $S$ and for all $x>0$:
$$\mathcal{A}_{\{\widehat{\mathcal{T}}\}}(S)= \mathcal{A}_{\{\omega_{\widehat{T}^x}\}}(S),$$
and the same equalities are valid for the corresponding sequence classes $\Lambda$.
So, $\mathcal{A}_{\{\widehat{\mathcal{T}}\}}(S)$ coincides with the space $\mathcal{A}_{\{\omega\}}(S)$ associated with a normalized weight $\omega=\omega_{\widehat{T}^x}$ satisfying
\hyperlink{om3}{$(\omega_3)$} and \hyperlink{om4}{$(\omega_4)$}; moreover, $\gamma(\omega)=\gamma(\tau)+1>1$.
\end{theorem}

\demo{Proof}
We do not wish to include here the details, since there is no significant difference with those carefully presented in \cite{testfunctioncharacterization} for a similar result in the ultradifferentiable case. The main idea behind the proof of this statement is that the ultraholomorphic classes considered here, associated either with a weight function or with a weight matrix, are introduced in exactly the same way as in the ultradifferentiable case, what lets us apply similar arguments as those developed in~\cite{dissertation,compositionpaper,testfunctioncharacterization} as long as only the structural properties of the spaces are concerned.

In order to be more specific, note that $\mathcal{T}$ satisfies \hyperlink{R-mg}{$(\mathcal{M}_{\{\on{mg}\}})$} and \hyperlink{R-L}{$(\mathcal{M}_{\{\on{L}\}})$} since it is associated with the weight $\tau$ (see \cite[Section 5]{compositionpaper}). Both properties are also true immediately for the weight matrix $\widehat{\mathcal{T}}$, and clearly each $\widehat{T}^x\in\hyperlink{LCset}{\mathcal{LC}}$, $x>0$.
Finally, by Corollary \ref{necessarytransformlemma1corollary}, we have every ingredient to mimic the proof of \cite[Corollary 3.17]{testfunctioncharacterization} in order to obtain the result,
except for the information about the index. For that, note first that Lemma~\ref{assofuncproperbis} implies that, for every $x>0$,
$\tau$ and $\omega_{T^x}$ are equivalent, so $\gamma(\omega_{T^x})=\gamma(\tau)>0$. Hence, by Corollary~\ref{necessarytransformlemma1corollary} and Lemma~\ref{lemmmaIndexLowerLegendreEnvelope} we have that
$$
\gamma(\omega_{\widehat{T}^x})=\gamma((\omega_{T^x}^{\iota})_{\star})= \gamma(\omega_{T^x})+1=\gamma(\tau)+1>1,
$$
as desired.
\qed\enddemo

\begin{remark}
The normalized weights satisfying \hyperlink{om3}{$(\omega_3)$}, \hyperlink{om4}{$(\omega_4)$} and with $\gamma(\omega)>1$ (this last condition amounts to \hyperlink{omsnq}{$(\omega_{\text{snq}})$}) are sometimes called \textit{strong weights}.
Observe that we also have $\mathcal{A}_{\{\widehat{\mathcal{T}}\}}(S)= \mathcal{A}_{\{(\omega_{T^x}^{\iota})_{\star}\}}(S)$ for any $x>0$,
but the weight function $(\omega_{T^x}^{\iota})_{\star}$ is concave and so it cannot satisfy the normalization condition; moreover, property \hyperlink{om4}{$(\omega_4)$} is also not clear for this weight.
\end{remark}

Until the end of this section we deal with the second situation, whose main aim is to prove the equivalence of the weight matrix associated with a strong weight and a weight matrix consisting of strongly log-convex sequences.

We start with the following result.

\begin{corollary}\label{coroWeightGammagreater1MatrixStronglyLogConvex}
Let $\omega$ be a normalized weight function satisfying \hyperlink{om3}{$(\omega_3)$}, \hyperlink{om4}{$(\omega_4)$} and $\gamma(\omega)>1$.
Let $\Omega=\{W^x:x>0\}$ be the weight matrix associated with $\omega$, and
define $L_p^x:=p!(w_p^x)^{\on{lc}}$. Then, for all $x>0$ it holds $\omega_{W^x}\hyperlink{sim}{\sim}\omega_{L^x}$.
\end{corollary}
\demo{Proof}
The condition $\gamma(\omega)>1$ guarantees that $\omega$ is equivalent to a concave weight function and also that $\omega$ has \hyperlink{om5}{$(\omega_5)$} (see Lemma~\ref{strongweightspace} and Proposition~\ref{Prop13MT88}).
Then, we are in a position to apply Lemma~\ref{lemmaWeightequivConcave} in order to conclude.
\qed\enddemo

\begin{corollary}\label{corocomparisonmaintheorem}
Let $\omega$ be a normalized weight function satisfying \hyperlink{om3}{$(\omega_3)$} and \hyperlink{om4}{$(\omega_4)$} with $\gamma(\omega)>1$, and let $\Omega=\{W^x: x>0\}$ be the associated weight matrix.  Consider the weight matrix $\widehat{\mathcal{T}}=\{\widehat{T}^x=(p!T^x_p)_{p\in\NN}: x>0\}$ with $T^x_p:=\exp(\frac{1}{x}\varphi^{*}_{\omega_{w^1}}(xp))$, $x>0$ and $p\in\NN$.
Then, for any sector $S$ we obtain (as locally convex vector spaces)
\begin{equation*}
\mathcal{A}_{\{\Omega\}}(S)=\mathcal{A}_{\{\omega\}}(S)=
\mathcal{A}_{\{\omega_{\widehat{T}^1}\}}(S)=
\mathcal{A}_{\{\widehat{\mathcal{T}}\}}(S),
\end{equation*}
and the same equalities are valid for the corresponding sequence classes $\Lambda$. Moreover, $\Omega$ is equivalent with respect to \hyperlink{Mroumapprox}{$\{\approx\}$} to the weight matrix $\widehat{\mathcal{T}}$ consisting only of strongly log-convex sequences.
\end{corollary}
\demo{Proof}
We put $\tau=\omega_{w^1}$ and observe that $T^1=(w^1)^{\on{lc}}$, so that $\widehat{T}_p^1=p!(w_p^1)^{\on{lc}}$. By Lemma~\ref{assofuncproperbis} and Corollary~\ref{coroWeightGammagreater1MatrixStronglyLogConvex} we have $\omega\hyperlink{sim}{\sim}\omega_{W^1}\hyperlink{sim}{\sim} \omega_{\widehat{T}^1}$,
and the argument is completed applying Theorem~\ref{theoEqualSpacesMatrixwithfactorialWeight} to the weight $\tau$ (the choice $\tau=\omega_{w^x}$ for any $x>0$ would also work).

\qed\enddemo

This corollary gives an answer to \cite[Question 5.11]{whitneyextensionmixedweightfunction}: For any given strong weight $\omega$ there does exist a weight matrix equivalent to $\Omega$ and satisfying \cite[$(1.3)$]{whitneyextensionmixedweightfunction} since for each strongly log-convex sequence $M$ the sequence of quotients $(\mu_p/p)_p$ is nondecreasing.\vspace{6pt}

In the next result we do not start with a weight function but with a sequence.
\begin{corollary}\label{mixedweightsequlemma}
Let $\widehat{M}=(p!M_p)_{p\in\NN}\in\hyperlink{LCset}{\mathcal{LC}}$ be given such that $\gamma(\omega_{\widehat{M}})>1$. $\mathcal{T}:=\{T^x: x>0\}$ shall denote the \hyperlink{Msc}{$(\mathcal{M}_{\on{sc}})$} matrix associated with the weight $\omega_M(=\omega_{M^{\on{lc}}})$, $\Omega:=\{W^x: x>0\}$ shall be the \hyperlink{Msc}{$(\mathcal{M}_{\on{sc}})$} matrix associated with the weight $\omega_{\widehat{M}}$ and finally consider $\widehat{\mathcal{T}}:=\{\widehat{T}^x: x>0\}$ defined by $\widehat{T}^x=(p!T^x_p)_{p\in\NN}$ for each $x>0$.


Then, for any sector $S$ we obtain (as locally convex vector spaces)
\begin{equation*}
\mathcal{A}_{\{\widehat{\mathcal{T}}\}}(S)=
\mathcal{A}_{\{\omega_{\widehat{T}^1}\}}(S)=
\mathcal{A}_{\{\omega_{\widehat{M}}\}}(S)=
\mathcal{A}_{\{\Omega\}}(S),
\end{equation*}
and the same equalities are valid for the corresponding sequence classes $\Lambda$.

Moreover $\widehat{\mathcal{T}}\hyperlink{Mroumapprox}{\{\approx\}}\Omega$ holds true, hence the matrix $\Omega$ is equivalent to a matrix consisting only of strongly log-convex weight sequences.
\end{corollary}

\demo{Proof}
It suffices to apply the previous result to $\omega=\omega_{\widehat{M}}$.
\qed\enddemo


In particular, Corollary \ref{mixedweightsequlemma} can be applied to each $\widehat{M}=(p!M_p)_{p\in\NN}\in\hyperlink{LCset}{\mathcal{LC}}$ satisfying \hyperlink{gamma1}{$(\gamma_1)$} (by recalling $(iv)$ in Lemma \ref{assofuncproper}).
%
%
But even for not strongly nonquasianalytic weight sequences we could apply Corollary \ref{mixedweightsequlemma} since, as mentioned above in Section \ref{assofunction}, it may happen that $\gamma(\omega_M)>1$ for a sequence $M$ not having \hyperlink{gamma1}{$(\gamma_1)$}.
We construct in \cite{JimenezSanzSchindlIndices} such an example explicitly.


\section{Existence of sectorially flat functions}\label{sectSectoriallyFlatFunct}

\subsection{Construction of outer functions}\label{constructionouterfunction}

The aim of this paragraph is to obtain holomorphic functions in the right half-plane whose growth is accurately controlled by a given weight function.

The next result transfers \hyperlink{omsnq}{$(\omega_{\text{snq}})$} for a weight function $\tau$ into a property for $\tau^{\iota}$, where $\tau^{\iota}(t)=\tau(1/t)$. Compare this with \cite[Lemma 2.1.1]{Thilliezdivision}.


\begin{lemma}\label{functionh3tau}
Let $\tau$ be a weight function. Then, one has $\gamma(\tau)>1$ if, and only if,
\begin{equation*}
\exists\;C\ge 1\;\forall\;y>0:\;\;\;\int_0^1 -\tau^{\iota}(ty)dt\ge-C(\tau^{\iota}(y)+1).
\end{equation*}
\end{lemma}

\demo{Proof}
First, by Lemma \ref{strongweightspace} we have $\gamma(\tau)>1$ if, and only if, the weight $\tau$ has \hyperlink{omsnq}{$(\omega_{\on{snq}})$}. In this condition we change $y\mapsto y^{-1}$, $t\mapsto t^{-1}$ and it is then equivalent to
\begin{equation}\label{functionh3equ1tau}
\exists\;C\ge 1\;\forall\;y>0:\;\;\;\int_0^1\tau\left(\frac{1}{ty}\right)dt\le C\tau\left(\frac{1}{y}\right)+C
\end{equation}
(observe that, by putting $s:=t^{-1}$, we get $\int_1^{\infty}\frac{\tau(ty)}{t^2}dt= \int_1^0\frac{\tau(s^{-1}y)}{s^{-2}}\left(-\frac{1}{s^2}ds\right) =\int_0^1\tau(s^{-1}y)ds$). Now we multiply \eqref{functionh3equ1tau} by $-1$ and recall that $\tau^{\iota}(t)=\tau(1/t)$.
\qed\enddemo

Moreover, we get the analogous result to \cite[Lemma 2.1.2]{Thilliezdivision}.

\begin{lemma}\label{Lemma212tau}
Let $\tau$ be a weight function such that $\gamma(\tau)>1$,
then
$$\int_{-\infty}^{+\infty}\frac{-\tau^{\iota}(|t|)}{1+t^2}dt>-\infty.$$
\end{lemma}

\demo{Proof}
Since by assumption Proposition \ref{Prop13MT88} and \cite[Corollary 1.4]{MeiseTaylor88} can be applied to $\tau$, there exists some $0<\alpha<1$ and $C\ge 1$ such that $\tau(t)\le Ct^{\alpha}+C$ for all $t>0$. Hence by definition and multiplying this by $-1$ we get $-\tau^{\iota}(t)\ge -C(t^{-\alpha}+1)$ for all $t>0$, from where the conclusion easily follows.
\qed\enddemo

In the next step we transfer \cite[Lemma 2.1.3]{Thilliezdivision} to the weight function case.

\begin{lemma}\label{Lemma213tau}
Let $\tau$ be a weight function
with $\gamma(\tau)>1$. Then for all $a>0$ there exists a function $F_a$ which is holomorphic on the right half-plane $H_1:=\{w\in\CC: \Re(w)>0\}$ and constants $A,B\ge 1$ (large) depending only on $\tau$ such that
\begin{equation}\label{Lemma213equtau}
\forall\;w\in H_1:\;\;\;B^{-a}\exp(-2a\tau^{\iota}(B^{-1}\Re(w)))\le|F_a(w)|\le \exp(-\frac{a}{2}\tau^{\iota}(A|w|)).
\end{equation}
\end{lemma}

\demo{Proof}
We are following the idea of the proof of \cite[Lemma 2.1.3]{Thilliezdivision}. For $w\in H_1$ put
$$F_a(w):=\exp\left(\frac{1}{\pi}\int_{-\infty}^{+\infty} \frac{-a\tau^{\iota}(|t|)}{1+t^2}\frac{itw-1}{it-w}dt\right);$$
Lemma \ref{Lemma212tau} implies immediately that $F_a$ is a holomorphic function in $H_1$. Since $F_a(w)=(F_1(w))^a$, we need only consider in the proof $a=1$ and put for simplicity $F:=F_1$.

For $w\in H_1$ write $w=u+iv$, hence $u>0$.
We have
$$\log(|F(w)|)=\frac{1}{\pi}\int_{\RR}-\tau^{\iota}(|t|) \frac{u}{(t-v)^2+u^2}dt=-\frac{1}{\pi}f\ast g_u(v),
$$
where $f(t):=\tau^{\iota}(|t|)$, $g_u(t):=u/(t^2+u^2)$. $f$ and $g_u$ are symmetrically nonincreasing functions, hence the convolution too. This means that $(f\ast g_u)(x)\le(f\ast g_u)(y)\le(f\ast g_u)(0)$ for $|x|\ge|y|\ge 0$. Consequently, the minimum for $w\mapsto\log(|F(w)|)$ is attained for $v=0$, so on the positive real axis and we have for all $w\in H_1$:
$$\log(|F(w)|)\ge\log(|F(u)|)=\log(|F(\Re(w))|), \hspace{10pt}\log(|F(u)|)=\frac{1}{\pi}\int_{\RR}-\tau^{\iota}(|t|) \frac{u}{t^2+u^2}dt=-\frac{1}{\pi}f\ast g_u(0).$$
First we concentrate on the left hand side in \eqref{Lemma213equtau}. Consider $K>0$ (small) and get
$$\pi\log(|F(u)|)=\int_{\RR}-\tau^{\iota}(|t|) \frac{u}{t^2+u^2}dt
=\int_{\{t:|t|\ge Ku\}}-\tau^{\iota}(|t|)\frac{u}{t^2+u^2}dt+\int_{\{t:|t|\le Ku\}}-\tau^{\iota}(|t|)\frac{u}{t^2+u^2}dt.$$
For the first integral we estimate by
$$\int_{\{|t|\ge Ku\}}-\tau^{\iota}(|t|)\frac{u}{t^2+u^2}dt
\ge-\tau^{\iota}(Ku)\int_{\{|t|\ge Ku\}}\frac{u}{t^2+u^2}dt
=-\tau^{\iota}(Ku)(\pi-2\arctan(K)),$$
since $t\mapsto -\tau^{\iota}(t)$ is nondecreasing.\vspace{6pt}

For the second integral we get
$$\int_{\{t: |t|\le Ku\}}-\tau^{\iota}(|t|)\frac{u}{t^2+u^2}dt
=\int_{\{s: |s|\le1\}}-\tau^{\iota}(Ku|s|)\frac{K}{K^2s^2+1}ds
\ge K\int_{\{s: |s|\le 1\}}-\tau^{\iota}(Ku|s|)ds,$$
since $-\tau^{\iota}(Ku|s|)\le 0$ holds for any $K,u>0$ and $|s|\le 1$. Let $C\ge 1$ be the constant appearing in Lemma \ref{functionh3tau}, then
$$K\int_{\{s:|s|\le 1\}}-\tau^{\iota}(Ku|s|)ds\ge2KC(-\tau^{\iota}(Ku)-1).$$
Thus we get for any $u>0$:
\begin{align*}
\pi\log(|F(u)|)&\ge(\pi-2\arctan(K)+2KC)(-\tau^{\iota}(Ku))-2KC
\ge(\pi+3(C-1)K)(-\tau^{\iota}(Ku))-2KC
\\&
\ge\pi(1+\frac{3}{\pi}(C-1)K)(-\tau^{\iota}(Ku))-2KC,
\end{align*}
since for all $K>0$ chosen sufficiently small enough but arbitrarily $\pi-2\arctan(K)+2KC$ behaves like $\pi+2(C-1)K+O(K^2)\le\pi+3(C-1)K$. Equivalently we have
$$\forall\;u>0:\;\;\;|F(u)|\ge
\exp(-2KC)(\exp(-\tau^{\iota}(Ku)))^{1+3(C-1)K/\pi}.$$
If $K>0$ is chosen small enough to have $1+\frac{3}{\pi}(C-1)K\le 2\Leftrightarrow K\le\frac{\pi}{3(C-1)}$, then since $\exp(-\tau^{\iota}(t))\le 1$ for any $t>0$ we get
\begin{align*}
\exp(-2KC)(\exp(-\tau^{\iota}(Ku)))^{1+3(C-1)K/\pi}&
\ge\exp(-2KC)(\exp(-\tau^{\iota}(Ku)))^2\\
&=
\exp(-2KC)\exp(-2\tau^{\iota}(Ku)).
\end{align*}
So the left hand side of \eqref{Lemma213equtau} is shown ($-\tau^{\iota}$ is nondecreasing).\vspace{6pt}

For the right hand side assume that $K>0$ is chosen arbitrarily (large), then
$$|F(w)|\!\le\!\exp\left(\frac{1}{\pi}\int_{\{t: |t-v|\le Ku\}}-\tau^{\iota}(|t|))\frac{u}{(t-v)^2+u^2}dt+\frac{1}{\pi}\int_{\{t: |t-v|\ge Ku\}}-\tau^{\iota}(|t|)\frac{u}{(t-v)^2+u^2}dt\right).$$
We estimate the second integral by $0$ (since the integrand is negative). For the first one, since we have $|t|\le|t-v|+|v|\le Ku+|w|\le K|w|+|w|=(1+K)|w|$ and $-\tau^{\iota}$ is nondecreasing, we get for any $w\in H_1$:
$$|F(w)|\le\exp\left(\frac{1}{\pi}(-\tau^{\iota}((1+K)|w|))\int_{\{t: |t-v|\le Ku\}}\frac{u}{(t-v)^2+u^2}dt\right).$$
Since the last integral is equal to $2\arctan(K)$ we summarize:
$$\forall\;w\in H_1:\;\;\;|F(w)|\le(\exp(-\tau^{\iota}((1+K)|w|)))^{2\arctan(K)/\pi}
\le(\exp(-\tau^{\iota}(((1+K)|w|)))^{1-2/(\pi K)}$$
holds, where $K>0$ is chosen sufficiently large to guarantee $\frac{2\arctan(K)}{\pi}\ge 1-\frac{2}{\pi K}$. If $K$ is chosen large enough to have $1-\frac{2}{\pi K}\ge\frac{1}{2}\Leftrightarrow K\ge\frac{4}{\pi}$, then we get
$$\forall\;w\in H_1:\;\;\;|F(w)|\le(\exp(-\tau^{\iota}(((1+K)|w|)))^{1/2}$$
which concludes the proof.
\qed\enddemo

\subsection{Construction of sectorially flat functions}\label{sectoriallyflatfunction}

Given a sequence $M\in\hyperlink{LCset}{\mathcal{LC}}$, it is easy to express flatness in the classes $\mathcal{A}_{\{M\}}(S)$ by means of the associated functions $\omega_M$ or $h_M$. Indeed, as in the classical Gevrey case, flat functions are characterized as those exponentially decreasing in a precise sense. The proof of the following result is a straightforward adaptation of the arguments in \cite[Proposition 4]{Thilliezsmoothsolution}.

\begin{lemma}\label{lemmaflatnessultraholclasses}
Let $S$ be an unbounded sector.
\begin{itemize}
\item[$(i)$] Let $M\in\mathcal{LC}$ be given such that $\lim_{p\to\infty}m_p^{1/p}=\infty$. Then,
\begin{itemize}
\item[$(i.1)$] If $f\in\mathcal{A}_{\{M\}}(S)$ is flat,
\begin{equation}\label{equaGlobalFlatnessSequence}
\exists C>0\;\exists k>0:\;\forall z\in S,\,\ |f(z)|\le Ch_m(k|z|)=C\exp(-\omega_{m}(1/(k|z|))).
\end{equation}
\item[$(i.2)$] Conversely, if $f$ is a holomorphic function in $S$ verifying~\eqref{equaGlobalFlatnessSequence}, then for every unbounded and proper subsector $T$ of $S$ one has
$f\in\mathcal{A}_{\{M\}}(T)$ and $f$ is flat.
\end{itemize}

\item[$(ii)$] Let $\mathcal{M}=\{M^x:x\in\RR_{>0}\}$ be a standard log-convex weight matrix with $\lim_{k\rightarrow\infty}(m^x_k)^{1/k}=+\infty$ for every $x>0$. Then,
\begin{itemize}
\item[$(ii.1)$] If $f\in\mathcal{A}_{\{\mathcal{M}\}}(S)$ is flat,
\begin{equation}\label{equaGlobalFlatnessMatrix}
\exists C>0\;\exists k>0\;\exists x>0:\;\forall z\in S,\ |f(z)|\le Ch_{m^x}(k|z|)=C\exp(-\omega_{m^x}(1/(k|z|))).
\end{equation}
\item[$(ii.2)$] Conversely, if $f$ is a holomorphic function in $S$ verifying~\eqref{equaGlobalFlatnessMatrix}, then for every unbounded and proper subsector $T$ of $S$ one has
$f\in\mathcal{A}_{\{\mathcal{M}\}}(T)$ and $f$ is flat.
\end{itemize}

\item[$(iii)$] Let $\omega$ be a normalized weight function with  $\hyperlink{om1}{(\omega_1)}$ and  $\hyperlink{om3}{(\omega_3)}$, and $\Omega=\{W^x=(W^x_j)_{j\in\NN}: x>0\}$ the associated weight matrix, for which (as indicated in Subsection~\ref{ultraholomorphicroumieu}) we have $\mathcal{A}_{\{\omega\}}(S)=\mathcal{A}_{\{\Omega\}}(S)$. Suppose moreover that $\lim_{p\rightarrow\infty}(w^x_p)^{1/p}=+\infty$ for every $x>0$. Then,
\begin{itemize}
\item[$(iii.1)$] If $f\in\mathcal{A}_{\{\omega\}}(S)$ is flat,
\begin{equation}\label{equaGlobalFlatnessMatrixfromWeightFunction}
\exists C>0\;\exists x>0:\;\forall z\in S,\ |f(z)|\le Ch_{w^x}(|z|)=C\exp(-\omega_{w^x}(1/|z|)).
\end{equation}
\item[$(iii.2)$] Conversely, if $f$ is a holomorphic function in $S$ verifying~\eqref{equaGlobalFlatnessMatrixfromWeightFunction}, then for every unbounded and proper subsector $T$ of $S$ one has
$f\in\mathcal{A}_{\{\omega\}}(T)$ and $f$ is flat.
\end{itemize}

\end{itemize}
\end{lemma}

\demo{Proof}
We only prove (i), since the rest of items may be obtained similarly. If $f\in\mathcal{A}_{\{M\}}(S)$,
there exists $k>0$ such that $f\in\mathcal{A}_{M,k}(S)$, and so for every $z\in S$ and $p\in\NN$ we have
\begin{equation}
   \label{equafbelongsAMhS}
   |f^{(p)}(z)|\le \Vert f\Vert_{M,k} k^pM_p.
\end{equation}
Now, by Taylor's formula we may write, for any $\lambda\in(0,1)$,
$$
f(z)-\sum_{j=0}^{p-1}\frac{f^{(j)}(\lambda z)}{j!}(z-\lambda z)^{j}=\frac{z^p}{(p-1)!}\int_{\lambda}^1 (1-t)^{p-1}f^{(p)}(tz)\,dt,
$$
and taking limits as $\lambda\to 0$ we deduce, since $f$ is flat, that
$$
f(z)=\frac{z^p}{(p-1)!}\int_{0}^1 (1-t)^{p-1}f^{(p)}(tz)\,dt.
$$
So, we may apply~\eqref{equafbelongsAMhS} in order to see that for every $p\in\NN$ we have
$$
|f(z)|\le \frac{|z|^p}{p!}\sup_{w\in S}|f^{(p)}(w)|\le
\Vert f\Vert_{M,k} (k|z|)^p m_p,
$$
what implies \eqref{equaGlobalFlatnessSequence} by definition of $h_m$ and~\eqref{functionhequ2}.

Conversely, suppose $f$ satisfies~\eqref{equaGlobalFlatnessSequence}. Given a proper and unbounded subsector $T$ of $S$, there exists $\varepsilon>0$ such that for every $z\in T$, the disc $D(z,\varepsilon|z|)$ is contained in $S$, and by Cauchy's formula we have, for every $p\in\NN$,
$$
f^{(p)}(z)=\frac{p!}{2\pi i}\int_{|w-z|=\varepsilon |z|}
\frac{f(w)}{(w-z)^{p+1}}\,dw,\quad p\in\NN.
$$
So, we easily estimate
$$
|f^{(p)}(z)|\le\frac{p!}{(\varepsilon |z|)^p}\max_{|w-z|=\varepsilon |z|}|f(w)|\le\frac{Cp!}{(\varepsilon |z|)^p}h_m(k(1+\varepsilon)|z|).
$$
Since $h_m(t)=\inf_{n\in\NN}m_nt^n$, on the one hand we have that
$$
|f^{(p)}(z)|\le\frac{Cp!}{(\varepsilon |z|)^p}m_p(k(1+\varepsilon)|z|)^p=
C\left(\frac{k(1+\varepsilon)}{\varepsilon}\right)^pM_p,
$$
from where $f\in\mathcal{A}_{M,k(1+\varepsilon)/\varepsilon}(T)\subset
\mathcal{A}_{\{M\}}(T)$, and on the other hand we deduce
$$
|f^{(p)}(z)|\le\frac{Cp!}{(\varepsilon |z|)^p}m_{p+1}(k(1+\varepsilon)|z|)^{p+1}=
C\frac{(k(1+\varepsilon))^{p+1}}{\varepsilon^p} \frac{M_{p+1}}{p+1} |z|,
$$
what immediately implies that $f^{(p)}(0)=\lim_{z\to 0,\,z\in T}f^{(p)}(z)=0$ for every $p\in\NN$, and $f$ is flat.
\qed\enddemo

\begin{remark}
The condition $\lim_{p\to\infty}m_p^{1/p}=\infty$ is not necessary for item (i) to hold. However, note that whenever $\lim_{p\to\infty}m_p^{1/p}<\infty$ the statement is trivial, since $h_m$ identically vanishes in an interval with 0 as its left-end point, and we immediately deduce that the only flat function in the class is the null function. Similar observations can be made for the other two items.
\end{remark}

\begin{remark}
Suppose given a normalized weight function $\omega$ with \hyperlink{om1}{$(\omega_1)$}, \hyperlink{om3}{$(\omega_3)$} and \hyperlink{om5}{$(\omega_5)$}, and let $\Omega=\{W^x=(W^x_p)_{p\in\NN}: x>0\}$ be its associated weight matrix. According to the information in Lemma~\ref{lemmaflatnessultraholclasses}(iii), one may characterize flatness in the ultraholomorphic class $\mathcal{A}_{\{\omega\}}(S)$ (which coincides with $\mathcal{A}_{\{\Omega\}}(S)$) in terms of exponential decrease of the type $\exp(-\omega_{w^x}(1/|z|))$ for some $x>0$. If $\omega$ has moreover \hyperlink{om4}{$(\omega_4)$}, by \eqref{omegaconjugateequivalentcorequ2} we see that this could be also expressed by exponential decrease of the type $\exp(-C\omega^{\star}(D|z|))$ for suitable $C,D>0$.
\end{remark}

Using the results from the previous sections the aim is now to transfer \cite[Theorem 2.3.1]{Thilliezdivision} to the weight function setting. Although the considered weights $\tau$ will satisfy \hyperlink{om1}{$(\omega_1)$}, we use the equivalent condition $\gamma(\tau)>0$ (see Lemma~\ref{lemmagammaomegapositiveimpliesomega1}), as this quantity will essentially indicate the opening of the sectors where our constructions will be valid.

\begin{theorem}\label{Theorem231tau}
Let $\tau$ be a weight function 
with $\gamma(\tau)>0$. Then for any $0<\gamma<\gamma(\tau)$ there exist constants $K_1,K_2,K_3>0$ depending only on $\tau$ and $\gamma$ such that for all $a>0$ there exists a function $G_a$ holomorphic in $S_{\gamma}$ and satisfying
\begin{equation}\label{Theorem231equtau}
\forall\;\xi\in S_{\gamma}:\;\;K_1^{-a}\exp(-2a\tau^{\iota}(K_2|\xi|)) \le|G_a(\xi)|\le
\exp(-\frac{a}{2}\tau^{\iota}(K_3|\xi|)).
\end{equation}
Moreover, if $\tau$ is normalized and satisfies  \hyperlink{om3}{$(\omega_3)$}, and $\mathcal{T}=\{T^x=(T^x_p)_{p\in\NN}:x>0\}$ is its associated weight matrix, then $G_a$ is a flat function in $\mathcal{A}_{\{\widehat{\mathcal{T}}\}}(S_{\gamma})$, where $\widehat{\mathcal{T}}$ is the standard log-convex weight matrix consisting of the sequences $\widehat{T}^x=(p!T^x_p)_{p\in\NN}$, $x>0$.

Finally, if we also assume that $\tau$ satisfies \hyperlink{om4}{$(\omega_4)$}, then there exist $x>0$ and $K_4>0$, both depending on $a$, such that
\begin{equation}\label{equaBoundsOptimalFlat}
\forall\;\xi\in S_{\gamma}: |G_a(\xi)|\ge K_4\,h_{T^{x}}(K_2|\xi|).
\end{equation}
\end{theorem}

We remark that \eqref{equaBoundsOptimalFlat} tells us that $G_a$ is indeed an optimal flat function, in the sense that its size is controlled by the functions $h_{T^x}$ not only from above, as needed for flatness, but also from below.

\demo{Proof} Let $a>0$ be arbitrary. Take $s,\delta>0$ such that $\gamma<\delta<\gamma(\tau)$, $s\delta<1<s\gamma(\tau)$.
By \eqref{newindex3} we get $s\gamma(\tau)>1\Leftrightarrow\gamma(\tau^{1/s})>1$, hence $\tau^{1/s}(t)=\tau(t^{1/s})=\tau^{\iota}(\frac{1}{t^{1/s}})
=(\tau^{\iota})^{1/s}(t)$ satisfies \hyperlink{omsnq}{$(\omega_{\on{snq}})$}. So we can use Lemma \ref{Lemma213tau} for the weight $\tau^{1/s}$ instead of $\tau$ and we obtain a function $F_a$ satisfying \eqref{Lemma213equtau} with  $\tau^{\iota}$ replaced by $(\tau^{\iota})^{1/s}$.
Then put
$$G_a(\xi)=F_a(\xi^s)\;\;\;\;\xi\in S_{\delta}.$$
Note that, as $s\delta<1$, the ramification $\xi\mapsto\xi^s$ maps holomorphically $S_{\delta}$ into $S_{\delta s}\subseteq S_1=H_1$, and so $G_a$ is well-defined.

We show that the restriction of $G_a$ to $S_{\gamma}\subseteq S_{\delta}$ satisfies the desired properties by proving that \eqref{Theorem231equtau} holds indeed on the whole $S_{\delta}$.

First we consider the lower estimate.
Let $\xi\in S_{\delta}$ be given, then $\Re(\xi^s)\ge\cos(s\delta\pi/2)|\xi|^s$ (since $s\delta\pi/2<\pi/2$). If $B\ge 1$ denotes the constant coming from the left hand side in \eqref{Lemma213equtau} applied for the weight $\tau^{1/s}$, then
\begin{align*}
|G_a(\xi)|&=|F_a(\xi^s)|\ge B^{-a}\exp(-2a(\tau^{\iota})^{1/s}(B^{-1}(\Re(\xi^s))))\\
&\ge B^{-a}\exp(-2a(\tau^{\iota})^{1/s}(B^{-1}\cos(s\delta\pi/2)|\xi|^s))
\\&
=B^{-a}\exp(-2a(\tau^{\iota})^{1/s}((B_1|\xi|)^s))
=B^{-a}\exp(-2a\tau^{\iota}(B_1|\xi|)),
\end{align*}
where we have put $B_1:=(B^{-1}\cos(s\delta\pi/2))^{1/s}$.\vspace{6pt}

Now we consider the right hand side in \eqref{Theorem231equtau} and proceed as before. Let $A$ be the constant coming from the right hand side of \eqref{Lemma213equtau} applied to $\tau^{1/s}$, so
\begin{align*}
|G_a(\xi)|&=|F_a(\xi^s)|\le \exp(-\frac{a}{2}(\tau^{\iota})^{1/s}(A|\xi|^s))
=\exp(-\frac{a}{2}(\tau^{\iota})^{1/s}((A^{1/s}|\xi|)^s))
\\&
=\exp(-\frac{a}{2}\tau^{\iota}(A^{1/s}|\xi|)),
\end{align*}
and \eqref{Theorem231equtau} has been proved for every $\xi\in S_{\delta}$.

Assume now that $\tau$ satisfies also \hyperlink{om3}{$(\omega_3)$}.
First put in the estimate above $A_1:=A^{1/s}$.
By using \eqref{equaExp-Omega.less.hm} for any $y>0$ we get
$$\exp(-\frac{a}{2}\tau^{\iota}(A_1|\xi|))\le h_{T^y}(A_1|\xi|)^{ya/2}.$$
Hence taking $y:=2a^{-1}$ proves that
\begin{equation}\label{equaGsubaflat}
\forall \xi\in S_{\delta}:\ \ |G_a(\xi)|\le h_{T^{2/a}}(A_1|\xi|),
\end{equation}
and it suffices to take into account Lemma~\ref{lemmaflatnessultraholclasses}.(ii.2) in order to deduce that $G_a$ belongs to $\mathcal{A}_{\{\widehat{\mathcal{T}}\}}(S_{\gamma})$ and it is flat.

Finally, if $\tau$ 
satisfies moreover \hyperlink{om4}{$(\omega_4)$} we may apply \eqref{equaExp-Omega.bigger.hm} for any $x>0$ and deduce that
$$\exp(-2a\tau^{\iota}(B_1|\xi|))\ge \exp(-2aC_x)h_{T^x}(B_1|\xi|)^{4xa}.$$
Now we take $x:=1/(4a)$ and prove that
$$
\forall \xi\in S_{\delta}:\ \ |G_a(\xi)|\ge B^{-a}\exp(-2aC_x)h_{T^{1/(4a)}}(B_1|\xi|),
$$
as desired.


\qed\enddemo

\section{Right inverses for the asymptotic Borel map in ultraholomorphic classes in sectors}\label{sectRightInver1var}

The aim of this section is to obtain an extension result in the ultraholomorphic classes considered. The existence of the flat functions $G_a$ obtained in Theorem~\ref{Theorem231tau} will be the main ingredient in the proof, which will follow the same technique as in previous works of A. Lastra, S. Malek and the second author~\cite{LastraMalekSanzContinuousRightLaplace, Sanzsummability}. Although for this construction the weight function $\tau$ needs not be normalized, we are interested in working with the weight matrix associated with it, which will be standard log-convex if we ask for normalization and  \hyperlink{om3}{$(\omega_3)$} to hold. Moreover, since the condition $\gamma(\tau)>0$ is also necessary and this amounts to \hyperlink{om1}{$(\omega_1)$}, we will have the warranty that the ultraholomorphic spaces associated with the weight function and its corresponding weight matrix coincide, see the comments preceding~\eqref{equaEqualitySpacesWeightFunctionMatrix}.

Note that any weight function may be substituted by a normalized equivalent one, and equivalence preserves the properties \hyperlink{om3}{$(\omega_3)$} and $\gamma(\tau)>0$, so it is no restriction to ask for normalization from the very beginning.

The next lemma provides us with suitable kernels entering the formal and analytic, truncated Laplace-like transforms we will need in our main statement.

\begin{lemma}\label{lemmaKernels}
Let $\tau$ be a normalized weight function with $\gamma(\tau)>0$ which satisfies
\hyperlink{om3}{$(\omega_3)$},  let $\mathcal{T}=\{T^x=(T^x_p)_{p\in\NN}:x>0\}$ be its associated weight matrix, let $0<\gamma<\gamma(\tau)$, and for $a>0$ let $G_a$ be the function constructed in Theorem~\ref{Theorem231tau}.
Let us define the function $e_{a}:S_{\gamma}\to\C$  by
$$
e_a(z):=z\,G_a(1/z),\quad z\in S_{\gamma}.
$$
The function $e_{a}$ enjoys the following properties:
\begin{itemize}
\item[(i)] $z^{-1}e_{a}(z)$ is uniformly integrable at the origin, it is to say, for any $t_0>0$ we have
    $$\sup_{|\sigma|<\gamma\pi/2} \int_{0}^{t_0}t^{-1}|e_{a}(te^{i\sigma})|dt<\infty.$$
\item[(ii)] There exist constants $K>0$, independent from $a$, and $C>0$, depending on $a$, such that
\begin{equation}\label{equaBounds.for.kernel}
|e_{a}(z)|\le Ch_{T^{4/a}}\left(\frac{K}{|z|}\right),\qquad z\in S_{\gamma}.
\end{equation}
\item[(iii)] For $\xi\in\R$, $\xi>0$, the values of $e_{a}(\xi)$ are positive real.
\end{itemize}
\end{lemma}


\demo{Proof}
(i) Let $t_0>0$ and $\sigma\in\R$ with $|\sigma|<\frac{\gamma\pi}{2}$. From \eqref{Theorem231equtau} we deduce that there exists $K_3>0$ such that
$$\int_{0}^{t_0}\frac{|e_{a}(te^{i\sigma})|}{t}dt\le \int_{0}^{t_0}\exp(-\frac{a}{2}\tau^{\iota}(K_3/t))dt\le t_0.$$

(ii) For the second part, we may apply
\eqref{equaGsubaflat} and write
$$|e_{a}(z)|=|z||G_{a}(1/z)|\le|z|h_{T^{2/a}}(A_1/|z|)$$
for every $z\in S_{\gamma}$, where $A_1$ does not depend on $a$.

We recall that from~\eqref{newmoderategrowth} we know that $2\omega_{T^{2x}}(t)\le \omega_{T^x}(t)$ for every $x>0$ and $t\ge 0$, and so $h_{T^x}(t)\le h_{T^{2x}}(t)^2$. Hence, combining this fact with the very definition of $h_{T^{2x}}$, we get
$$
|e_{a}(z)|\le |z|h_{T^{4/a}}(A_1/|z|)^2\le |z|\Big(\frac{A_1}{|z|}\Big)T_1^{4/a}h_{T^{4/a}}(A_1/|z|)
< A_1T_1^{4/a}h_{T^{4/a}}(A_1/|z|),
$$
as desired.

(iii) Finally, if $\xi>0$ then $e_{a}(\xi)=\xi G_{a}(1/\xi)$. From the integral expression for $G_a$, 
it is immediate to check that the imaginary part of the integrand is an odd function, so the imaginary part of $G_a(1/\xi)$ is 0, while the real part is positive.
\qed\enddemo


\begin{definition}
We define the \textit{moment function} associated with the function $e_a$ (introduced in the previous Lemma) as
$$m_a(\lambda):=\int_{0}^{\infty}t^{\lambda-1}e_{a}(t)dt=
\int_{0}^{\infty}t^{\lambda}G_{a}(1/t)dt.$$
\end{definition}

From Lemma~\ref{lemmaKernels} and the definition of $h_{T^x}$ we see that for every $p\in\NN$,
$$
|e_{a}(z)|\le C\frac{K^pT^{4/a}_p}{|z|^p},\qquad z\in S_{\gamma}.
$$
So, we easily deduce that the function $m_a$ is well defined and continuous in $\{\hbox{Re}(\lambda)\ge0\}$, and holomorphic in $\{\hbox{Re}(\lambda)>0\}$. Moreover, $m_a(\xi)$ is positive  for every $\xi\ge0$, and the sequence $(m_a(p))_{p\in\NN}$ is called the \textit{sequence of moments} of $e_{a}$.

The next result is similar to Proposition 3.6 in~\cite{LastraMalekSanzContinuousRightLaplace}. The fact that such estimates could also be obtained in the present situation became clear thanks to the arguments by O. Blasco in~\cite{Blasco}.

\begin{proposition}\label{propmequivm}
Let $\tau$ be a normalized weight function with $\gamma(\tau)>0$ which satisfies 
\hyperlink{om3}{$(\omega_3)$}, let $\mathcal{T}=\{T^x=(T^x_p)_{p\in\NN}:x>0\}$ be its associated weight matrix, and for $0<\gamma<\gamma(\tau)$ and $a>0$ let $G_a,e_a,m_a$ be the functions previously constructed. Then, there exist constants $C_1,C_2>0$, both depending on $a$, such that for every $p\in\NN$ one has
\begin{equation}
  \label{equaIneqMomentsMatrix}
  C_1(\frac{K_2}{2})^pT_p^{1/(2a)}\le m_a(p)\le C_2K_3^pT_{p}^{4/a},
\end{equation}
where $K_2$ and $K_3$ are the constants, not depending on $a$, appearing in Theorem~\ref{Theorem231tau}.
\end{proposition}
\demo{Proof}
Let $p\in\NN_0$. By the second inequality in~\eqref{Theorem231equtau}, we have for every $s>0$ that
$$m_a(p)=\int_{0}^{\infty}t^{p}G_a(1/t)dt\le \int_{0}^{s}t^{p}dt
+\int_{s}^{\infty}\frac{1}{t^2}\exp\big((p+2)\log(t) -\frac{a}{2}\tau(\frac{t}{K_3})\big)dt.
$$
Now, observe that
\begin{align*}
\sup_{t>s}\big((p+2)\log(t)-\frac{a}{2}\tau(\frac{t}{K_3})\big)
&= (p+2)\log(K_3)+ \sup_{u>s/K_3}\big((p+2)\log(u)-\frac{a}{2}\tau(u)\big)\\
&\le (p+2)\log(K_3)+ \sup_{u>0}\big((p+2)\log(u)-\frac{a}{2}\tau(u)\big)\\
&=(p+2)\log(K_3)+ \frac{a}{2}\sup_{v\in\RR}\big(\frac{2(p+2)}{a}v-\tau(e^v)\big)
\\&=(p+2)\log(K_3)+ \frac{a}{2}\varphi_{\tau}^{*}(\frac{2(p+2)}{a}).
\end{align*}
Hence, we deduce that
$$m_a(p)\le\frac{s^{p+1}}{p+1}
+K_3^{p+2} \exp(\frac{a}{2}\varphi_{\tau}^{*}(\frac{2(p+2)}{a})) \frac{1}{s}.
$$
Since this is valid for any $s>0$, we compute the infimum of such bounds as $s$ runs in $(0,\infty)$, whose value is
$$
\frac{p+2}{p+1}K_3^{p+1}
\exp(\frac{a}{2}\frac{p+1}{p+2} \varphi_{\tau}^{*}(\frac{2(p+2)}{a})),
$$
and obtain that
$$m_a(p)\le 2K_3^{p+1}\exp(\frac{a}{2} \varphi_{\tau}^{*}(\frac{2(p+2)}{a}))=
2K_3^{p+1}T_{p+2}^{2/a}\le
(2K_3T_2^{4/a})K_3^pT_p^{4/a},
$$
where we have made use of~\eqref{newmoderategrowth}.

For the second part of the estimates we use the first inequality in~\eqref{Theorem231equtau} and the fact that $\tau$ is nondecreasing in order to write, for every $s>0$,
$$m_a(p)\ge \int_{0}^{s}t^{p}G_a(1/t)dt\ge K_1^{-a}\int_{0}^{s}t^{p}\exp\big( -2a\tau(\frac{t}{K_2})\big)dt\ge
K_1^{-a}\frac{s^{p+1}}{p+1}\exp\big( -2a\tau(\frac{s}{K_2})\big).
$$
We compute now the supremum of such bounds and, in a similar way, we deduce that
\begin{align*}
m_a(p)&\ge \frac{K_1^{-a}}{p+1}\exp\sup_{s>0} \big((p+1)\log(s)-2a\tau(\frac{s}{K_2})\big)=
\frac{K_1^{-a}}{p+1}K_2^{p+1}\exp(2a \varphi_{\tau}^{*}(\frac{p+1}{2a}))\\
&=\frac{K_1^{-a}}{p+1}K_2^{p+1}T_{p+1}^{1/(2a)}\ge
(K_1^{-a}K_2)(\frac{K_2}{2})^p T_p^{1/(2a)},
\end{align*}
as desired.

\qed\enddemo


The proof of the incoming result  rests on a constructive procedure which combines a formal Borel transform and a truncated Laplace transform, like the original one in the Gevrey case (see~\cite{Tougeron},~\cite[Theorem\ 4.1]{SanzGevreyseveralvar}). The main tool needed is a suitable kernel, namely the function $e_a$ obtained in Lemma~\ref{lemmaKernels}, in terms of which both aforementioned transforms are explicitly given. This generalizes the classical situation, where the role of $e_a$ is played by the exponential function, whose moment function is precisely the Euler Gamma function.

\begin{theorem}\label{theoExtensionOperatorsMatrix}
Let $\tau$ be a normalized weight function with $\gamma(\tau)>0$ which satisfies 
\hyperlink{om3}{$(\omega_3)$}, let $0<\gamma<\gamma(\tau)$, let $\mathcal{T}=\{T^x=(T^x_p)_{p\in\NN}:x>0\}$ be its associated weight matrix, and consider the weight matrix $\widehat{\mathcal{T}}=\{\widehat{T}^x:x>0\}$ where $\widehat{T}^x=(p!T^x_p)_{p\in\NN}$. Then, there exists a constant $k_0>0$ such that for every $x>0$ and every $h>0$, one can construct a linear and continuous map
$$
\lambda\in\Lambda_{\widehat{T}^x,h}\mapsto f_{\lambda}\in \mathcal{A}_{\widehat{T}^{8x},k_0h}(S_{\gamma})
$$
such that for every $\lambda$ one has $\mathcal{B}(f_{\lambda})=\lambda$.

In particular, the Borel map
$\mathcal{B}:\mathcal{A}_{\{\widehat{\mathcal{T}}\}}(S_{\gamma})\to \Lambda_{\{\widehat{\mathcal{T}}\}}$ is surjective.
\end{theorem}

\demo{Proof}
Fix $\delta>0$ such that $\gamma<\delta<\gamma(\tau)$. Given $\lambda=(\lambda_p)_{p\in\NN}\in\Lambda_{\widehat{T}^x,h}$, we have
\begin{equation}\label{equaNormlambda}
|\lambda_p|\le |\lambda|_{\widehat{T}^x,h} h^{p}p!T_p^x,\quad p\in\N_0.
\end{equation}
We choose $a=1/(2x)$, and consider the function $G_a$, defined in $S_{\delta}$, obtained in Theorem~\ref{Theorem231tau} for such value of $a$, and the corresponding functions $e_a$ and $m_a$ previously defined.
Next, we consider the formal power series
$$
\widehat{f}_{\lambda}:=\sum_{p=0}^\infty\frac{\lambda_p}{p!}z^p
$$
and its formal (Borel-like) transform
$$
\widehat{\mathcal{B}}_a\widehat{f}_{\lambda}:= \sum_{p=0}^\infty\frac{\lambda_p}{p!m_a(p)}z^p.
$$
By the choice of $a$, \eqref{equaNormlambda} and the first part of the inequalities in~\eqref{equaIneqMomentsMatrix}, we deduce that
\begin{equation}
  \label{equaSizeCoefficientsBorelTransform}
  \left|\frac{\lambda_p}{p!m_a(p)}\right|\le
  \frac{|\lambda|_{\widehat{T}^x,h} h^{p}p!T_p^x}
  {C_1(K_2/2)^p p!T_p^x}=
  \frac{|\lambda|_{\widehat{T}^x,h}}{C_1}
  \Big(\frac{2h}{K_2}\Big)^p,
\end{equation}
and so the series $\widehat{\mathcal{B}}_a\widehat{f}_{\lambda}$ converges in the disc of center 0 and radius $K_2/(2h)$ (not depending on $\lambda$), where it defines a holomorphic function $g_\lambda$. We set $R_0:=K_2/(4h)$, and define
\begin{equation*}
f_{\lambda}(z):=\int_{0}^{R_0}e_{a}\left(\frac{u}{z}\right) g_{\lambda}(u)\frac{du}{u},\qquad z\in S_{\delta}.
\end{equation*}
By virtue of Leibniz's theorem on analyticity of parametric integrals, $f_\lambda$ is holomorphic in $S_{\delta}$.

Our next aim is to obtain suitable estimates for the difference between $f$ and the partial sums of the series $\widehat{f}_{\lambda}$.

Let $N\in\N_0$ and $z\in S_{\delta}$. We have
\begin{align*}
f_{\lambda}(z)-\sum_{p=0}^{N-1}\lambda_p\frac{z^p}{p!} &= f_{\lambda}(z)-\sum_{p=0}^{N-1}\frac{\lambda_p}{m_a(p)}m_a(p)\frac{z^p}{p!}\\
&= \int_{0}^{R_0}e_{a}\left(\frac{u}{z}\right) \sum_{p=0}^{\infty}\frac{\lambda_{p}}{m_a(p)}\frac{u^p}{p!} \frac{du}{u} -\sum_{p=0}^{N-1}\frac{\lambda_p}{m_a(p)} \int_{0}^{\infty}u^{p-1}e_{a}(u)du\frac{z^p}{p!}.
\end{align*}
In the second integral we make the change of variable $v=zu$, what results in a rotation of the line of integration. By the estimate~(\ref{equaBounds.for.kernel}), one may use Cauchy's residue theorem in order to obtain that
$$
z^p\int_{0}^{\infty}u^{p-1}e_{a}(u)du= \int_{0}^{\infty}v^{p}e_{a}\left(\frac{v}{z}\right)\frac{dv}{v},
$$
which allows us to write the preceding difference as
\begin{multline*}
\int_{0}^{R_0}e_{a}\left(\frac{u}{z}\right) \sum_{p=0}^{\infty}\frac{\lambda_{p}}{m_a(p)} \frac{u^p}{p!}\frac{du}{u} -\sum_{p=0}^{N-1}\frac{\lambda_p}{m_a(p)} \int_{0}^{\infty}u^{p}e_{a}\left(\frac{u}{z}\right)\frac{du}{u} \frac{1}{p!}\\
=\int_{0}^{R_0}e_{a}\left(\frac{u}{z}\right) \sum_{p=N}^{\infty}\frac{\lambda_{p}}{m_a(p)}\frac{u^p}{p!} \frac{du}{u} -\int_{R_0}^{\infty}e_{a}\left(\frac{u}{z}\right) \sum_{p=0}^{N-1}\frac{\lambda_p}{m_a(p)} \frac{u^{p}}{p!}\frac{du}{u}.
\end{multline*}
Then, we have
\begin{equation}\label{equaBoundsAsympExpanflambda}
  \left|f_{\lambda}(z)-\sum_{p=0}^{N-1}\lambda_p\frac{z^p}{p!}\right|\le |f_{1}(z)|+|f_2(z)|,
\end{equation}
where
$$f_{1}(z)=\int_{0}^{R_0}e_{a}\left(\frac{u}{z}\right) \sum_{p=N}^{\infty}\frac{\lambda_{p}}{m_a(p)} \frac{u^p}{p!}\frac{du}{u},\quad
f_{2}(z)=\int_{R_0}^{\infty}e_{a}\left(\frac{u}{z}\right) \sum_{p=0}^{N-1}\frac{\lambda_p}{m_a(p)} \frac{u^{p}}{p!}\frac{du}{u}.$$
From~\eqref{equaSizeCoefficientsBorelTransform} we deduce that
\begin{align}
|f_{1}(z)|&\le \frac{|\lambda|_{\widehat{T}^x,h}}{C_1}
  \int_{0}^{R_0}\left|e_{a}\left(\frac{u}{z}\right)\right| \sum_{p=N}^{\infty}\Big(\frac{2hu}{K_2}\Big)^p\frac{du}{u}
  =\frac{|\lambda|_{\widehat{T}^x,h}}{C_1} \Big(\frac{2h}{K_2}\Big)^N
  \int_{0}^{R_0}\left|e_{a}\left(\frac{u}{z}\right)\right| \frac{u^N}{1-\frac{2hu}{K_2}}\frac{du}{u}\nonumber\\
  &\le \frac{2|\lambda|_{\widehat{T}^x,h}}{C_1} \Big(\frac{2h}{K_2}\Big)^N
  \int_{0}^{R_0}\left|e_{a}\left(\frac{u}{z}\right)\right| u^{N-1}\,du,\label{equaBoundsf-1}
\end{align}
where in the last step we have used that $0<u<R_0=K_2/(4h)$ we have $1-2hu/K_2>1/2$.
In order to estimate $f_{2}(z)$, observe that for $u\ge R_0$ and $0\le p\le N-1$ we always have $u^p\le R_0^pu^N/R_0^N$, and so, using again~\eqref{equaSizeCoefficientsBorelTransform} and the value of $R_0$, we may write
$$\left|\sum_{p=0}^{N-1}\frac{\lambda_pu^p}{p!m_a(p)}\right|\le \frac{|\lambda|_{\widehat{T}^x,h}}{C_1}\frac{u^N}{R_0^N} \sum_{p=0}^{N-1}\Big(\frac{2h}{K_2}\Big)^p R_0^p
\le \frac{2|\lambda|_{\widehat{T}^x,h}}{C_1} \Big(\frac{4h}{K_2}\Big)^N u^N.
$$
Then, we deduce that
\begin{equation}\label{equaBoundsf-2}
  |f_2(z)|\le \frac{2|\lambda|_{\widehat{T}^x,h}}{C_1} \Big(\frac{4h}{K_2}\Big)^N \int_{R_0}^{\infty}\left|e_{a}\left(\frac{u}{z}\right)\right| u^{N-1}du.
\end{equation}
In order to conclude, it suffices then to obtain estimates for $\int_{0}^{\infty}|e_{a}(u/z)|u^{N-1}du$. For this, note first that, by the estimates in~\eqref{Theorem231equtau},
\begin{align*}
  \int_{0}^{\infty}\left|e_{a}\left(\frac{u}{z}\right)\right| u^{N-1}du&=
  \int_{0}^{\infty}\frac{u}{|z|} \left|G_{a}\left(\frac{z}{u}\right)\right| u^{N-1}du\\
  &\le \int_0^\infty \frac{u^N}{|z|}\exp\big(-\frac{a}{2}\tau(\frac{u}{K_3|z|}) \big)\,du=
  |z|^N\int_0^\infty t^N\exp\big(-\frac{a}{2}\tau(\frac{t}{K_3})\big)\,dt.
\end{align*}
Now, we can follow the first part of the proof of Proposition~\ref{propmequivm} to obtain that
\begin{equation}
  \label{equaBoundsIntegralsKerneltimespower}
  \int_{0}^{\infty}\left|e_{a}\left(\frac{u}{z}\right)\right| u^{N-1}du\le C_2K_3^NT_N^{4/a}|z|^N=C_2K_3^NT_N^{8x}|z|^N.
\end{equation}
Gathering~\eqref{equaBoundsAsympExpanflambda}, \eqref{equaBoundsf-1}, \eqref{equaBoundsf-2} and \eqref{equaBoundsIntegralsKerneltimespower}, we get
\begin{equation}
  \label{equaAsExflambdafinal}
  \left|f_{\lambda}(z)-\sum_{p=0}^{N-1}\lambda_p\frac{z^p}{p!}\right|\le \frac{2C_2|\lambda|_{\widehat{T}^x,h}}{C_1} \Big(\frac{4hK_3}{K_2}\Big)^NT_N^{8x}|z|^N.
\end{equation}
A straightforward application of Cauchy's integral formula for the derivatives (as in the proof of Lemma~\ref{lemmaflatnessultraholclasses}) shows that there exists a constant $r$, depending only on $\gamma$ and $\delta$, such that whenever $z$ is restricted to belong to $S_{\gamma}$, one has that
for every $p\in\NN$,
\begin{equation*}
|f^{(p)}(z)|\le \frac{2C_2|\lambda|_{\widehat{T}^x,h}}{C_1} \Big(\frac{4hK_3r}{K_2}\Big)^p p!T_p^{8x}.
\end{equation*}
So, putting $k_0:=\frac{4K_3r}{K_2}$ (independent from $x$ and $h$), we see that
$f_{\lambda}\in \mathcal{A}_{\widehat{T}^{8x},k_0h}(S_{\gamma})$ and $\Vert f_{\lambda}\Vert_{\widehat{T}^{8x},k_0h}\le
\frac{2C_2}{C_1}|\lambda|_{\widehat{T}^x,h}$. Since the map sending $\lambda$ to $f_{\lambda}$ is clearly linear, this last inequality implies that the map is also continuous from $\Lambda_{\widehat{T}^x,h}$ into
$\mathcal{A}_{\widehat{T}^{8x},k_0h}(S_{\gamma})$. Finally, from \eqref{equaAsExflambdafinal} one may easily deduce that
$\mathcal{B}(f_{\lambda})=\lambda$, and we conclude.
\qed\enddemo

\begin{remark}
Indeed, the estimates in \eqref{equaAsExflambdafinal} show precisely that the function $f_\lambda$ admits the series $\widehat{f}_\lambda$ as its uniform asymptotic expansion in the sector $S_\delta$, with constraints given mainly in terms of the sequence $T^{8x}$. The link between the classes of functions admitting such an expansion and the ultraholomorphic classes studied in this paper is extremely strong, as it can be seen in \cite{Sanzflatultraholomorphic}.
\end{remark}

We may infer also the existence of extension operators in the classes associated with the weight functions corresponding to the weight matrices $\widehat{\mathcal{T}}$.

\begin{corollary}\label{coroExtensionOperatorsWeights}
Let $\tau$ be a normalized weight function with $\gamma(\tau)>0$ which satisfies 
\hyperlink{om3}{$(\omega_3)$} and \hyperlink{om4}{$(\omega_4)$}, let $\gamma$, $\mathcal{T}$ and $\widehat{\mathcal{T}}$ be as in the previous Theorem, and let $\omega$ 
be the weight function given in Theorem~\ref{theoEqualSpacesMatrixwithfactorialWeight}, in such a way that $\mathcal{A}_{\{\widehat{\mathcal{T}}\}}(S_\gamma)= \mathcal{A}_{\{\omega\}}(S_\gamma)$.
Then, for every $l>0$ there exists $l_1>0$ such that there exists a linear and continuous map
$$
\lambda\in\Lambda_{\omega,l}\mapsto f_{\lambda}\in \mathcal{A}_{\omega,l_1}(S_{\gamma})
$$
such that for every $\lambda$ one has $\mathcal{B}(f_{\lambda})=\lambda$.
\end{corollary}

\demo{Proof}
Let $\Omega:=\{W^x: x>0\}$ be the weight matrix associated with the weight function $\omega$, i.e. $W^x_p:=\exp\left(\frac{1}{x}\varphi^{*}_{\omega}(xp)\right)$, and $\widehat{\mathcal{T}}:=\{\widehat{T}^x: x>0\}$, where $\widehat{T}^x_p:=p!T^x_p$ for each $x>0$ and $p\in\NN$.
We may apply \eqref{equaEqualitySpacesWeightFunctionMatrix} in order to deduce that
$\mathcal{A}_{\{\widehat{\mathcal{T}}\}}(S_\gamma)= \mathcal{A}_{\{\Omega\}}(S_\gamma)$. It turns out that, independently and by similar arguments, related only to the way the classes are defined, the same equality will hold for the corresponding ultradifferentiable spaces, introduced in \cite[Chapter\ 7]{dissertation} (see also \cite[4.2]{compositionpaper}). Now, Theorem 4.6 in \cite{compositionpaper} states that this equality in the ultradifferentiable case amounts to the equivalence of the corresponding weight matrices, in the sense that
\begin{align}
\label{equaFirstpartEquivalenceMatrices}
\forall\;x>0&\;\exists\;y>0\;\exists\;C>0:\forall p\in\NN\;\;\;W_p^x\le C^p p!T_p^y,\text{ and}\\
\label{equaSecondpartEquivalenceMatrices}\forall\;y>0&\;\exists\;x>0\;\exists\;D>0:\forall p\in\NN\;\;\;p!T_p^y\le D^p W_p^x.
\end{align}
We fix $l>0$. By \eqref{equaFirstpartEquivalenceMatrices}, there exist $x>0$ and $C_1>0$ such that for every $p\in\NN$ one has $W_p^l\le C_1^p p! T_p^x$. So, given $\lambda\in\Lambda_{\omega,l}$, for every $p\in\NN$ we have
$$
|\lambda_p|\le |\lambda|_{\omega,l}W_p^l\le |\lambda|_{\omega,l}C_1^p p! T_p^x,
$$
what implies that $\lambda\in\Lambda_{\widehat{T}^x,C_1}$ and $|\lambda|_{\widehat{T}^x,C_1}\le |\lambda|_{\omega,l}$. Now, consider the function $f_{\lambda}$ given by the previous theorem, which belongs to $\mathcal{A}_{\widehat{T}^{8x},k_0C_1}(S_{\gamma})$, depends on $\lambda$ in a linear continuous way (so, there exists $A>0$ with $\Vert f_{\lambda}\Vert_{\widehat{T}^{8x},k_0C_1}\le A|\lambda|_{\widehat{T}^x,C_1}$), and is such that $\mathcal{B}(f_{\lambda})=\lambda$. By \eqref{equaSecondpartEquivalenceMatrices} there exists $l_0>0$ and $C_2>0$ (independent from $\lambda$) such that for every $p\in\NN$, $\widehat{T}^{8x}_p\le C_2^pW_p^{l_0}$,
and by property \hyperlink{R-L}{$(\mathcal{M}_{\{\text{L}\}})$} for $\Omega$, there exist $l_1>0$ and $D>0$ such that $(k_0C_1C_2)^pW_p^{l_0}\le DW_p^{l_1}$. Then, we obtain that for every $p\in\NN$ and every $z\in S_{\gamma}$,
\begin{align*}
|f^{(p)}(z)|&\le \Vert f_{\lambda}\Vert_{\widehat{T}^{8x},k_0C_1}(k_0C_1)^p\widehat{T}_p^{8x}\le
\Vert f_{\lambda}\Vert_{\widehat{T}^{8x},k_0C_1}(k_0C_1C_2)^pW_p^{l_0}\\
&\le D\Vert f_{\lambda}\Vert_{\widehat{T}^{8x},k_0C_1}W_p^{l_1}\le
AD|\lambda|_{\widehat{T}^x,C_1}W_p^{l_1}\le AD|\lambda|_{\omega,l}W_p^{l_1}.
\end{align*}
This means that $f\in\mathcal{A}_{\omega,l_1}(S_{\gamma})$ and, moreover,
$\Vert f_{\lambda}\Vert_{\omega,l_1}\le AD|\lambda|_{\omega,l}$, so that the map is linear and continuous between the corresponding spaces.
\qed\enddemo

Using Corollary~\ref{corocomparisonmaintheorem} we can prove now our second main extension theorem, here we are starting with the weight matrix ``containing factorials''.

\begin{theorem}\label{comparisonmaintheorem}
Let $\omega$ be a normalized weight function satisfying \hyperlink{om3}{$(\omega_3)$}, \hyperlink{om4}{$(\omega_4)$} and $\gamma(\omega)>1$, and let $\Omega=\{W^x: x>0\}$ be the associated weight matrix. Then for all $0<\gamma<\gamma(\omega)-1$ and for all $x>0$ there exist $y>0$ and a continuous linear extension operator
\begin{equation}\label{comparisonthmequ0}
E^{\omega}_{\gamma,x}:\Lambda_{\omega,x}\longrightarrow\mathcal{A}_{\omega,y}(S_{\gamma}),
\end{equation}
i.e. $(\mathcal{B}E^{\omega}_{\gamma,x})(\lambda)=\lambda$ holds true for any sequence $\lambda\in\Lambda_{\omega,x}$.
In particular, the map
$\mathcal{B}:\mathcal{A}_{\{\omega\}}(S_{\gamma})\longrightarrow\Lambda_{\{\omega\}}$ is surjective.
\end{theorem}

\demo{Proof}
First, by Corollary \ref{growthindexgammacorollary}, we obtain $\gamma((\omega^{\star})^{\iota})=\gamma(\omega_{w^x})
=\gamma(\omega_{W^x})-1=\gamma(\omega)-1$, hence $\gamma(\omega_{w^x})>0$, $\omega_{w^x}$ has \hyperlink{om1}{$(\omega_1)$} and moreover $(\omega^{\star})^{\iota}\hyperlink{sim}{\sim}\omega_{w^x}$ for all $x>0$.
Thus, each weight $\omega_{w^x}$ satisfies the requirements to apply the first main extension result, Theorem \ref{theoExtensionOperatorsMatrix}, and we do it for the choice $\tau:=\omega_{w^1}$.
We get that for given $0<\gamma<\gamma(\tau)=\gamma(\omega)-1$ there does exist some $k_0>0$ such that for any $h>0$ and $x>0$ we can construct a linear and continuous map
\begin{equation}\label{comparisonthmequ2}
E_{\gamma,h}:\Lambda_{\widehat{T}^x,h}\longrightarrow\mathcal{A}_{\widehat{T}^{8x},k_0h}(S_{\gamma}),\hspace{30pt}\lambda\mapsto f_{\lambda},
\end{equation}
such that $\mathcal{B}(f_{\lambda})=\lambda$ for all $\lambda\in\Lambda_{\widehat{T}^x,h}$.
Let $l>0$ be given, arbitrary but from now on fixed and let $\lambda\in\Lambda_{\omega,l}\equiv\Lambda_{W^l,1}$. We follow the proof of Corollary \ref{coroExtensionOperatorsWeights} and use the equivalence $\widehat{\mathcal{T}}\hyperlink{Mroumapprox}{\{\approx\}}\Omega$ to get $W^l_j\le h^j\widehat{T}^x_j$ for some $x,h>0$ and all $j\in\NN$ resp. $(k_0h)^j\widehat{T}^{8x}_j\le(k_0hh_1)^jW^{l_1}_j\le CW^{l_2}_j$ for some $h_1,C>0$ and indices $l_1,l_2>0$ and all $j\in\NN$, where in the last step we have used \eqref{newexpabsorb}. Hence, using \eqref{comparisonthmequ2}, we have obtained now an extension operator
\begin{equation*}
E^{\omega}_{\gamma,l}:\Lambda_{\omega,l}\longrightarrow\mathcal{A}_{\omega,l_2}(S_{\gamma}),
\end{equation*}
i.e. $(\mathcal{B}E^{\omega}_{\gamma,l})(\lambda)=\lambda$ holds true for any sequence $\lambda\in\Lambda_{\omega,l}$ as desired.
\qed\enddemo

\begin{remark}
Unfortunately, and as it already happened in Corollary~\ref{coroExtensionOperatorsWeights},
we do not know whether in \eqref{comparisonthmequ0} we can have $y=cx$ for some $c>0$ not depending on given $x$. More precisely, the equivalence $\widehat{\mathcal{T}}\hyperlink{Mroumapprox}{\{\approx\}}\Omega$ does not imply necessarily $x=cl$, $l_1=c_18x$, moreover $h$, $h_1$ both are depending on given $l$, hence $l_2=Al_1$ for some $A$ also depending on $l$ (see \eqref{newexpabsorb}).
\end{remark}

\subsection{Application to a mixed setting}\label{subsectMixedSetting}

As commented in the introduction, the known extension results for Denjoy-Carleman ultraholomorphic classes of Roumieu type in unbounded sectors by V. Thilliez~\cite{Thilliezdivision} or J. Schmets and M. Valdivia~\cite{Schmetsvaldivia00} impose growth conditions on the weight sequence defining the classes, namely moderate growth in the first case, and $(\beta_2)$ condition (see~\eqref{equaConditionBeta2}) in the second one.
The aim in this last paragraph is to indicate how our previous results may be used in order to obtain extension results in a mixed setting under minimal assumptions on the sequence. We will discuss two situations:

\begin{itemize}
\item[(a)] Let us consider a weight sequence $\widehat{M}$ which is \hyperlink{lc}{\text{(lc)}} and has $\hyperlink{gamma1}{(\gamma_1)}$.
As a consequence of the results by H.-J. Petzsche (see~\cite{petzsche}), $\widehat{M}$ may be substituted by a strongly equivalent sequence $\widehat{L}$
which is \hyperlink{slc}{\text{(slc)}} and also has $\hyperlink{gamma1}{(\gamma_1)}$. We write $\widehat{L}=(n!L_n)_{n\in\N_0}$, in such a way that $L:=(L_n)_{n\in\NN}\in\hyperlink{LCset}{\mathcal{LC}}$ and $\gamma(L)=\gamma(\widehat{L})-1=\gamma(\widehat{M})-1>0$ (see Subsection~\ref{growthindexgamma}).\par

Consider now the associated weight function for $L$, $\tau:=\omega_{L}$, which is a normalized weight function with \hyperlink{om3}{$(\omega_3)$} and \hyperlink{om4}{$(\omega_4)$}, and satisfies $\gamma(\tau)\ge\gamma(L)>0$ (see Remark~\ref{remarkComparisonIndicesGamma}). Let $\mathcal{T}=\{T^x=(T^x_p)_{p\in\NN}:x>0\}$ be its associated weight matrix, and consider the weight matrix $\widehat{\mathcal{T}}=\{\widehat{T}^x:x>0\}$ where $\widehat{T}^x=(p!T^x_p)_{p\in\NN}$. It turns out that $T^1=L$, and so $\widehat{T}^1=\widehat{L}$.\par
Hence, by Theorem~\ref{theoExtensionOperatorsMatrix} for any $0<\gamma<\gamma(\tau)$
there exists $k_0>0$ such that for every $h>0$, one can construct a linear and continuous map
\begin{equation}\label{equaExtenOperMixedSetting}
\lambda\in\Lambda_{\widehat{L},h}\equiv\Lambda_{\widehat{T}^1,h}\mapsto f_{\lambda}\in \mathcal{A}_{\widehat{T}^{8},k_0h}(S_{\gamma})
\end{equation}
such that for every $\lambda$ one has $\mathcal{\mathcal{B}}(f_{\lambda})=\lambda$.
Moreover, since $\widehat{M}\hyperlink{simeq}{\simeq}\widehat{L}$ we have $\Lambda_{\{\widehat{M}\}}=\Lambda_{\{\widehat{L}\}}$, and we deduce that $\mathcal{\mathcal{B}}(\mathcal{A}_{\{\widehat{T}^{8}\}}(S_{\gamma}))\supset \Lambda_{\{\widehat{M}\}}$.

\item[(b)] Suppose now that $\widehat{M}\in\hyperlink{LCset}{\mathcal{LC}}$ and $\gamma(\omega_{\widehat{M}})>1$, then we may apply Theorem~\ref{comparisonmaintheorem}. An example of this situation is presented in~\cite{JimenezSanzSchindlIndices} in which $\gamma(\widehat{M})=1$, so $\widehat{M}$ does not have $\hyperlink{gamma1}{(\gamma_1)}$, what excludes the possibility of applying another extension results.
\end{itemize}

\begin{remark}\label{remarkResultThilliezwithmg}
It is interesting to note that, in case the previously considered weight sequence $\widehat{M}$ has \hyperlink{lc}{\text{(lc)}}, $\hyperlink{gamma1}{(\gamma_1)}$ and \hyperlink{mg}{(\text{mg})}, we recover exactly the extension result of V. Thilliez~\cite[Theorem 3.2.1]{Thilliezdivision}. To see this, first note that \hyperlink{mg}{(\text{mg})} is stable under (weak or strong) equivalence of sequences, so the sequence $\widehat{L}$ in the previous item (a)
will also have \hyperlink{mg}{(\text{mg})}, and the same will hold for the sequence $L$, as indicated in Subsection~\ref{subsectionWeightSequences}.
So, on one hand, by Lemma~\ref{assofuncproper}.(iii) the weight function $\tau$ has \hyperlink{om6}{$(\omega_6)$}, and Remark~\ref{importantremark} implies that the matrix $\mathcal{T}$, and consequently also the matrix $\widehat{\mathcal{T}}$, is constant, in the sense that all the weight sequences appearing in it are equivalent to each other. This means then that $\widehat{M}$ is equivalent to $\widehat{T}^{8}$. On the other hand, as indicated in Remark~\ref{remarkComparisonIndicesGamma} and again by the moderate growth of $L$, we have $\gamma(\tau)=\gamma(L)=\gamma(M)$, and so the previous extension operator in~\eqref{equaExtenOperMixedSetting}, existing whenever $0<\gamma<\gamma(M)=\gamma(\widehat{M})-1$, can be seen as
$$
\lambda\in\Lambda_{\widehat{M},h}\mapsto f_{\lambda}\in \mathcal{A}_{\widehat{M},k_1h}(S_{\gamma})
$$
for some suitable $k_1>0$ not depending on $h$. In this case, $\mathcal{B}:\mathcal{A}_{\{\widehat{M}\}}(S_{\gamma})\to \Lambda_{\{\widehat{M}\}}$ is
surjective.
This is precisely the form of the extensions provided in Theorem 3.2.1 in~\cite{Thilliezdivision}.
\end{remark}

The next example shows that there do exist sequences for which previously known extension results by V. Thilliez or J. Schmets and M. Valdivia cannot be applied.

\begin{example}
\label{examSequenceBadGrowth}
We first recall that condition \hyperlink{gamma1}{$(\gamma_1)$} for a log-convex weight sequence $M$ is equivalent to the following condition (see \cite[Proposition 1.1]{petzsche}):
$$\hypertarget{beta1}{(\beta_1)}:\Leftrightarrow\exists\;k\in\NN_{>1}:
\liminf_{p\rightarrow\infty}\frac{\mu_{kp}}{\mu_p}>k.$$

Also in \cite{petzsche} the following condition was introduced:
\begin{equation}\label{equaConditionBeta2}
\hypertarget{beta2}{(\beta_2)}:\Leftrightarrow\;\forall\;
\varepsilon>0\;\exists\;k\in\NN_{>1}:\;
\limsup_{p\rightarrow\infty}\left(\frac{M_{kp}}{M_p}\right)^{\frac{1}{p(k-1)}}
\frac{1}{\mu_{kp}}\le\varepsilon.
\end{equation}
In \cite[Lemma 2.2 $(b)$]{Schmetsvaldivia00} it was pointed out that, by Stirling's formula, $M$ has \hyperlink{beta2}{$(\beta_2)$} if, and only if, $m$ has \hyperlink{beta2}{$(\beta_2)$}.

We show now that there exist sequences $M\in\RR_{>0}^{\NN}$ which satisfy $m\in\hyperlink{LCset}{\mathcal{LC}}$ and such that:
\begin{itemize}
\item[$(i)$] \hyperlink{beta2}{$(\beta_2)$},  \hyperlink{beta1}{$(\beta_1)$} and \hyperlink{mg}{(\text{mg})} are violated.
\item[$(ii)$] \hyperlink{beta2}{$(\beta_2)$} and \hyperlink{mg}{(\text{mg})} are violated, 
    but nevertheless \hyperlink{beta1}{$(\beta_1)$} holds.
\end{itemize}

$(i)$
We define $m:=(m_p)_p$ by putting $m_p:=q^{f(p)}$, where $q\ge\exp(1)$ and $f:[0,+\infty)\rightarrow[0,+\infty)$ is a convex function with $f(0)=0$ and $\lim_{p\rightarrow\infty}f(p)=+\infty$ defined as follows:

Let $(a_j)_{j\ge 1}$ be an increasing sequence in $\NN_{>0}$ with $a_{j+1}\ge a_j\cdot j$ for all $j\in\NN_{>0}$. Furthermore denote by $G^s=(G^s_p)_p$ the Gevrey-sequence, i.e. $G^s_p=p!^s$, $s>1$. Consider now the set of points
$$\mathcal{P}:=\{(a_j,j\log(a_j!))\}.$$
For $j\ge 1$ let $L_j$ be the line connecting the points $(a_j,j\log(a_j!))$ and $(a_{j+1},(j+1)\log(a_{j+1}!))$ with slope $l_j:=\frac{(j+1)\log(a_{j+1}!)-j\log(a_j!)}{a_{j+1}-a_j}$. For $j=0$ let $L_0$ be the line connecting $(0,0)$ with the point $(a_1,\log(a_1!))$.

By the log-convexity of $G^1$, the points on the line $L_0$ lie above each point on $\{(p,\log(p!)): 0<p<a_1\}$. By choosing $(a_j)_j$ increasing fast enough we can achieve that $(l_j)_j$ is increasing:

For this note that
$l_j\ge\tilde{l}_j$, where $\tilde{l}_j$ is the slope of the straight line $\tilde{L}_j$ connecting the points $(a_j,j\log(a_j!))$ and $(a_{j+1},j\log(a_{j+1}!))$. By convexity and the properties of the Gevrey-sequences $\tilde{l}_j$ is tending to infinity as $a_{j+1}\rightarrow\infty$, and so one recursively choose the $a_j$ in such a way that $l_{j+1}\ge \tilde{l}_{j+1}\ge l_j$.\vspace{6pt}

We put $f(p)$ equal to the height of the segment $L_j$ at the point $p$ for all $p\in\NN$ with $a_j\le p\le a_{j+1}$.

Since $(l_j)_j$ is increasing, $m$ is log-convex (and so $M$ is strongly log-convex). Moreover by construction and convexity $f(p)\ge j\log(p!)$ for all $p\ge a_j$, which proves $m_p\ge G^j_p$ for all $p\ge a_j$ (and arbitrary $j\in\NN_{>0}$). Since
$$
\forall\;
h>0\;\exists\;C\ge 1\;\forall\;p\in\NN:\; G^j_p\le Ch^pG^{j+1}_p,
$$
we get that
$$\forall\;
h>0\;\exists\;C\ge 1\;\forall\;p\in\NN:\; G^s_p\le Ch^pm_p
$$
for any $s>1$. So the sequence $m$ is ``beyond all Gevrey sequences'', and this fact excludes \hyperlink{mg}{(\text{mg})} (see~\cite{matsumoto}) for $m$ and, equivalently, for $M$.

Let us see that condition \hyperlink{beta2}{$(\beta_2)$} does not hold:

The expression in this condition gives $q^{S(k,p)}$ with $S(k,p):=\frac{1}{p(k-1)}(f(kp)-f(p))-(f(kp)-f(kp-1))$. So, in fact $S(k,p)$ is measuring the difference of two different slopes of $\{(p,f(p)): p\in\NN\}$. For any $k\in\NN$ and for all $j\ge k$ we get $ka_j\le ja_j\le a_{j+1}$, so $S(k,a_j)=0$, what implies that
$$\limsup_{p\rightarrow\infty}q^{S(k,p)}\ge
\limsup_{j\rightarrow\infty}q^{S(k,a_j)}=q^0=1,$$
what excludes \hyperlink{beta2}{$(\beta_2)$}.
Analogously we see that \hyperlink{beta1}{$(\beta_1)$} does not hold either: for any $p$ such that $a_j\le p-1<p\le a_{j+1}$ we have that $m_p/m_{p-1}=e^{l_j}$, and this implies that for any given $k\in\NN$, and whenever $j\ge k$, we have $ka_j\le ja_j\le a_{j+1}$ and $$\frac{\mu_{ka_j}}{\mu_{a_j}}=
\frac{ka_jm_{ka_j}/m_{ka_j-1}}{a_jm_{a_j}/m_{a_j-1}}=
k,
$$
whence $\liminf_{p\rightarrow\infty}\frac{\mu_{kp}}{\mu_p}=k$.
So, $\gamma(m)=0$ follows in this case.

It is worthy to comment that another example in this situation is the sequence $\widehat{M}$ mentioned in the previous item (b) before Remark~\ref{remarkResultThilliezwithmg}, and which is included in~\cite{JimenezSanzSchindlIndices}. $\widehat{M}$ is \hyperlink{slc}{\text{(slc)}} and does not have $\hyperlink{gamma1}{(\gamma_1)}$, so that it does not have $\hyperlink{beta2}{(\beta_2)}$ either (see~\cite[p.\ 223]{Schmetsvaldivia00}). Moreover, $\widehat{M}$ does not have \hyperlink{mg}{(\text{mg})}.

$(ii)$ However,  from the previously constructed sequence $M$ (or the sequence $\widehat{M}$ in (b)) it is possible to get a sequence with \hyperlink{beta1}{$(\beta_1)$}, and without \hyperlink{beta2}{$(\beta_2)$} and \hyperlink{mg}{(\text{mg})}.

For this we point out the following: Let $M$ be log-convex, then $\liminf_{p\rightarrow\infty}\frac{\mu_{kp}}{\mu_p}\ge 1$ holds.
Then,
the sequence $P:=(p!^2M_p)_p$ always satisfies \hyperlink{beta1}{$(\beta_1)$}: Note $\pi_p=p^2\mu_p$, hence $\liminf_{p\rightarrow\infty}\frac{\pi_{kp}}{\pi_p}=
k^2\cdot\liminf_{p\rightarrow\infty}\frac{\mu_{kp}}{\mu_p}\ge k^2>k$.
On the other hand, if $M$ is the sequence in (i), $M$ does not satisfy \hyperlink{beta2}{$(\beta_2)$} either \hyperlink{mg}{(\text{mg})}, and the same is true for $P$ since, as already commented, these two properties are stable under multiplication by the factorials.

\end{example}

\textbf{Acknowledgements}: The first two authors are partially supported by the Spanish Ministry of Economy, Industry and Competitiveness under the project MTM2016-77642-C2-1-P. The first author is partially supported by the University of Valladolid through a Predoctoral Fellowship (2013 call) co-sponsored by the Banco de Santander. The third author is supported by FWF-Project J~3948-N35, as a part of which he is an external researcher at the Universidad de Valladolid (Spain) for the period October 2016-September 2018.\par
The authors wish to express their gratitude to Prof. \'Oscar Blasco, from the Universidad de Valencia (Spain), for his helpful comments regarding Proposition~\ref{propmequivm}.

\bibliographystyle{plain}
\bibliography{Bibliography}

\vskip1cm

\textbf{Affiliation}:\\
J.~Jim\'{e}nez-Garrido, J.~Sanz:\\
Departamento de \'Algebra, An\'alisis Matem\'atico, Geometr{\'\i}a y Topolog{\'\i}a, Universidad de Valladolid\\
Facultad de Ciencias, Paseo de Bel\'en 7, 47011 Valladolid, Spain.\\
Instituto de Investigaci\'on en Matem\'aticas IMUVA\\
E-mails: jjjimenez@am.uva.es (J.~Jim\'{e}nez-Garrido), jsanzg@am.uva.es (J. Sanz).
\\
\vskip.5cm
G.~Schindl:\\
Departamento de \'Algebra, An\'alisis Matem\'atico, Geometr{\'\i}a y Topolog{\'\i}a, Universidad de Valladolid\\
Facultad de Ciencias, Paseo de Bel\'en 7, 47011 Valladolid, Spain.\\
E-mail: gerhard.schindl@univie.ac.at.

\end{document}